\documentclass[12pt]{article}
\usepackage{amsfonts}
\usepackage{amssymb}
\usepackage{latexsym}
\usepackage{amsmath}
\usepackage{bibmods}
\usepackage{graphicx,epsfig}
\usepackage{color}
\usepackage[arrow,matrix,curve]{xy}
\usepackage{amssymb}
\usepackage{amscd}
\usepackage{color}

\newtheorem{lemma}{\indent Lemma}[section]
\newtheorem{theorem}[lemma]{\indent Theorem}

\newtheorem{proposition}[lemma]{\indent Proposition}

\newtheorem{definition}[lemma]{\indent Definition}

\newtheorem{remark}[lemma]{\indent Remark}

\newcommand{\bt}{\beta}

\newcommand{\vf}{\varphi}
\newcommand{\ve}{\varepsilon}

\newcommand{\bB}{\mathbb{B}}
\newcommand{\bC}{\mathbb{C}}
\newcommand{\bH}{\mathbb{H}}
\newcommand{\bL}{\mathbb{L}}

\newcommand{\bR}{\mathbb{R}}
\newcommand{\bS}{\mathbb{S}}
\newcommand{\bW}{\mathbb{W}}
\newcommand{\bX}{\mathbb{X}}

\newcommand{\cB}{{\mathcal B}}
\newcommand{\cC}{{\mathcal C}}
\newcommand{\cD}{{\mathcal D}}
\newcommand{\cF}{{\mathcal F}}

\newcommand{\cK}{{\mathcal K}}

\newcommand{\cM}{{\mathcal M}}

\newcommand{\cS}{{\mathcal S}}

\newcommand{\cx}{{}{\scriptstyle{\mathcal X}}} 

\newcommand{\fB}{\mathfrak{B}}
\newcommand{\fM}{\mathfrak{M}}

\newcommand{\dst}{\displaystyle}
\newcommand{\pa}{\partial}
\newcommand{\ov}{\overline}
\newcommand{\wt}{\widetilde}
\newcommand{\dR}{\overset{\text{\tiny$\bullet$}}{\mathbb{R}}}

\newcommand{\A}{{\boldsymbol A}}

\newcommand{\D}{{\boldsymbol D}}
\newcommand{\F}{{\boldsymbol F}}

\newcommand{\K}{{\boldsymbol K}}
\newcommand{\M}{{\boldsymbol M}}
\newcommand{\N}{{\boldsymbol N}}
\newcommand{\R}{{\boldsymbol R}}

\newcommand{\V}{{\boldsymbol V}}
\newcommand{\W}{{\boldsymbol W}}

\newcommand{\nub}{{\boldsymbol{\nu}}}

\newcommand{\supp}{\operatorname{supp}}
\newcommand{\QED}{\hspace{\fill}$\Box$\medskip\par}

 \title{\Large MIXED BOUNDARY VALUE PROBLEMS FOR THE LAPLACE-BELTRAMI EQUATION}

\author{R. Duduchava  \& M. Tsaava}

\begin{document}
 \date{}
 \maketitle

\begin{abstract}
We investigate the mixed Dirichlet-Neumann boundary value problems for the Laplace-Beltrami equation on a smooth bounded surface $\mathcal{C}$ with a smooth boundary in non-classical setting in the Bessel potential space $\mathbb{H}^s_p(\mathcal{C})$ for $s>\dst\frac1p$, $1<p<\infty$. To the initial BVP we apply a quasi-localization and obtain a model BVP for the Laplacian. The model mixed BVP on the half plane is reduced to an equivalent system of Mellin convolution equation (MCE) in Sobolev-Slobode\v{c}kii space (potential method). MCE is ivestigated in both Bessel potential and Sobolev-Slobode\v{c}kii spaces. The symbol of the obtained system is written explicitly and is responsible for the Fredholm properties and the index of the system. An explicit criterion for the unique solvability of the initial BVP in the non-classical setting is derived as well.
\end{abstract}



\noindent {\bf ﻿Funding:} The research was supported by Shota Rustaveli National Science Foundation grants no. 13/14 and 31/39.

\vskip-5mm
\tableofcontents

\section*{Introduction and formulation of the problems}
\label{section0}
\setcounter{equation}{0}

Let $\cS\subset\bR^3$ be some closed orientable surface, bordering a compact inner $\Omega^+$ and outer $\Omega^-:=\bR^3\setminus\ov{\Omega^+}$ domains. By $\cC$ we denote a subsurface of $\cS$, which has two faces $\cC^-$ and $\cC^+$ and inherits the orientation from $\cS$: $\cC^+$ borders the inner domain $\Omega^+$  and $\cC^-$ borders the outer domain $\Omega^-$.  $\cC$ has the smooth boundary $\Gamma:=\pa\cC$, which is decomposed into two closed parts $\Gamma=\Gamma_D\cup\Gamma_N$, consisting each of finite number of smooth arcs, having in common only endpoints.

Let $\nub (\omega)=\left(\nu_1(\omega),\nu_2(\omega), \nu_3(\omega)\right)^\top$, $\omega\in\ov\cC$ be the unit normal vector field on the surface $\cC$ and $\partial_\nub=\dst\sum_{j=1}^3\nu_j\pa_j$ be the normal derivative. Let us consider the Laplace-Beltrami operator in $\mathcal{C}$ written in terms of the G\"unter's tangent derivatives (see \cite{DMM06,Du09,DTT14} for more details)
\begin{eqnarray}\label{e0.0}
\Delta_\cC:=\cD^2_1+\cD^2_2+\cD^2_3,\qquad \cD_j:=\pa_j-\nu_j\pa_\nub,\quad j=1,2,3.
\end{eqnarray}
Let $\nub_\Gamma(t)=({\nu_{\Gamma,1}(t),\nu_{\Gamma,2}(t), \nu_{\Gamma,3}(t)})^\top$, $t\in\Gamma$, be the unit normal vector field on the boundary $\Gamma$, which is tangential to the surface $\cC$ and directed outside of the surface. And, finally, let $\partial_{\nub_\Gamma}:=\dst\sum_{j=1}^3\nu_{\Gamma,j}\cD_j$ be the normal derivative on the boundary of the surface, which is the outer tangential derivative on the surface.

We study the following mixed boundary value problem for the Laplace-Beltrami equation
\begin{eqnarray}\label{e0.1}
\left\{\begin{array}{ll}
\Delta_\cC u(t)=f(t),\qquad & t\in\cC, \\[0.2cm]
u^+(\tau)=g(\tau),     \qquad & \tau\in\Gamma_D, \\[0.2cm]
(\partial_{\nub_\Gamma}u)^+(\tau)=h(\tau),\qquad & \tau\in\Gamma_N.
\end{array}\right.
\end{eqnarray}
where $u^+$ and $(\partial_{\nub_\Gamma}u)^+$ denote respectively the Dirichlet and the Neumann traces on the boundary.

We need the Bessel potential $\bH^s_p(\cS)$,  $\bH^s_p(\cC)$, $\wt{\bH}^s_p(\cC)$ and Sobolev-Slobode\v{c}kii $\bW^r_p(\Gamma)$ spaces, where $\cS$ is a closed smooth surface (without boundary), which contains $\cC$ as a subsurface, $1<p<\infty,\quad \frac1p<s<1-\frac1p$. The Bessel potential space $\bH^s_p(\bR^n)$ is defined as a subset of the space of Schwartz distributions $\bS'(\bR^n)$ endowedp with the norm (see \cite{Tr95})
 \[
  ||u\big|\bH_p^s(\bR^n)||:=||\langle D\rangle^su\big| L_p(\bR^n)||,
 \]
where $\langle D\rangle^s:=\cF^{-1}(1+|\xi|^2)^{\frac s2}\cF$ is the Bessel potential and $\cF$, $\cF^{-1}$ are the Fourier transformations.  For the definition of the Sobolev-Slobode\v{c}kii space $\bW_p^s(\bR^n)=\bB_{p,p}^s(\bR^n)$ see \cite{Tr95}. The space $\bW_p^s(\cS)$ coincides with the trace space of $\bH_p^{s+\frac1p}(\bR^3)$ on $\cS$ and is known that $\bW^s(\cS)=\bH^s(\cS)$ for $s\geq0$, $1<p<\infty$ (see \cite{Tr95}).

We use, as common, the notation $\bH^s(\cS)$ and $\bW^s(\cS)$ for the spaces $\bH_2^s(\cS)$ and $\bW_2^s(\cS)$ (the case $p=2$).

The spaces $\bH_p^s(\cS)$ and $\bW_p^s(\cS)$ are defined by a partition of the unity $\{\psi_j\}_{j=1}^\ell$ subordinated to some covering $\{ Y_j\}_{j=1}^\ell$ of $\cS$ and local coordinate diffeomorphisms (see \cite{Tr95,HW08} for details)
\[
\varkappa_j : X_j\rightarrow Y_j , \qquad  X_j\subset\bR^2\, ,\quad  j=1,\ldots,\ell.
\]

The space $\wt {\bH}_p^s(\cC)$ is defined as the subspace of $\bH_p^s(\cS)$ of those functions $\vf\in \bH_p^s(\cS)$, which are supported in the closed sub-surface $\supp\vf\subset\ov{\cC}$, whereas $\bH_p^s(\cC)$ denotes the quotient space $\bH_p^s(\cC):=\bH_p^s(\cS)\Big/\wt{\bH}_p^s(\cC^c)$, and $\cC^c:=\cC\setminus\ov{\cC}$ is the complemented sub-surface. For $s>1/p-1$ the space $\bH_p^s(\cC)$ can be identified with the space of those distributions $\vf$ on $\bR^n_+$ which admit extensions $\ell\vf\in\bH_p^s(\cS)$, while $\bH_p^s(\cC)$ is identified with the space $r_\cC\bH_p^s(\cS)$, where $r_\cC$ denotes the restriction from $\cS$ to  the sub-surface $\cC$.

It is worth noting that for an integer $m=1,2,\ldots$ the Sobolev spaces $\bH^m_p(\cS)$ and $\bW^m_p(\cS)$ coincide and the equivalent norm is defined with the help of the G\"unter's derivatives (see \cite{Du01,Du09,DMM06}):
 \[
  ||u\big|\bW_p^m(\cS)||:=\left[\sum_{\alpha|\leqslant m}||\cD^\alpha u\big| L_p(\cS)||^p\right]^{\frac1p},
       \quad\mbox{ where }\quad \cD^\alpha:=\cD^{\alpha_1}_1\cD^{\alpha_2}_2\cD^{\alpha_3}_3
 \]
and the G\"unter's derivatives $\cD_1, \cD_2, \cD_3$ are defined in \eqref{e0.0}.

 Let us also consider $\widetilde{\mathbb{H}}^{-1}_0(\mathcal{C})$, a subspace of $\widetilde{\mathbb{H}}^{-1}(\mathcal{C})$, orthogonal to
 \[
 \widetilde{\bH}^{-1}_\Gamma (\cC):=\left\{f \in\widetilde{\mathbb{H}}^{-1} (\cC)\;:\;\langle f,\varphi\rangle = 0 \;\text{for all}\;\varphi\in C^1_0(\cC)\right\}.
 \]
 $ \widetilde{\mathbb{H}}^{-1}_\Gamma(\mathcal{C})$ consists of those distributions on $\cS$, belonging to  $\widetilde{\mathbb{H}}^{-1}(\mathcal{C})$ which have their supports just on $\Gamma$ and  $\widetilde{\mathbb{H}}^{-1}(\mathcal{C})$ can be decomposed into the direct sum of subspaces:
  \[
 \widetilde{\mathbb{H}}^{-1}(\mathcal{C})= \widetilde{\mathbb{H}}^{-1}_\Gamma(\mathcal{C})\oplus
       \widetilde{\mathbb{H}}^{-1}_0(\mathcal{C}).
  \]
  The space $\widetilde{\bH}^{-1}_\Gamma (\cC)$ is non-empty (see \cite[\S\, 5.1]{HW08}) and excluding it from $\widetilde{\mathbb{H}}^{-1}(\mathcal{C})$ is needed to make  BVPs uniquelly solvable (cf. \cite{HW08} and the next Theorem \ref{t0.1}).

The Lax-Milgram Lemma applied to the BVP \eqref{e0.1} gives the following result.

\begin{theorem}[Theorem 14, \cite{DTT14} and \S\, 5.1, \cite{HW08}]\label{t0.1}
The BVP \eqref{e0.1} has a unique solution in the classical weak setting:
\begin{equation}\label{e0.2}
u\in\mathbb{H}^1(\mathcal{C}),\quad f\in\widetilde{\mathbb{H}}^{-1}_0(\mathcal{C}),
     \quad g\in\mathbb{H}^{1/2}(\Gamma_D), \quad h\in\mathbb{H}^{-1/2}(\Gamma_N).
\end{equation}
\end{theorem}

From Theorem \ref{t0.1} we can not even conclude that a solution is continuous. If we can prove that there is a solution $u\in\mathbb{H}^1_p(\mathcal{C})$ for some $2<p<\infty$, we can enjoy even a H\"older continuity of $u$. It is very important to know maximal smoothness of a solution as, for example, while designing approximation methods. To this end we will investigate the solvability properties of the BVP \eqref{e0.1} in the following non-classical setting
\begin{eqnarray}\label{e0.3}
u\in\mathbb{H}^s_p(\mathcal{C}),\quad
     f\in\widetilde{\mathbb{H}}^{s-2}_p(\mathcal{C})\cap\widetilde{\mathbb{H}}^{-1}_0(\mathcal{C}),\quad g\in\mathbb{W}^{s-1/p}_p(\Gamma),\quad h\in\mathbb{W}^{s-1-1/p}_p(\Gamma),\\
               1<p<\infty, \quad s>\frac1p\nonumber
\end{eqnarray}
and find necessary and sufficient conditions of solvability. Note, that the constraint $s>\dst\frac1p$ is necessary to ensure the existence of the trace $u^+$ on the boundary.

To formulate the main theorem of the present work we need the following definition.
 %
\begin{definition}\label{t0.2}
The BVP \eqref{e0.1}, \eqref{e0.3} is Fredholm if the homogeneous problem $f=g=h=0$ has a finite number of linearly independent solutions and only a finite number of orthogonality conditions on the data $f,g,h$ ensure the solvability of the BVP.
\end{definition}

We prove below the following theorem (see the concluding part of \S\ 5).
 %
\begin{theorem}\label{t0.3}
Let $1<p<\infty$, $s>\dst\frac1p$.

The BVP \eqref{e0.1} is Fredholm in the non-classical setting \eqref{e0.3} if and only if:
 \begin{eqnarray}\label{e0.4}
p\not=2\quad or \quad p=2\quad and \quad  s\not=\frac12+k, \qquad \text{for}\quad  k=0,1,2,\ldots.
 \end{eqnarray}
In particular, the BVP \eqref{e0.1} has a unique solution $u$ in the non-classical setting \eqref{e0.3} if
 \begin{eqnarray}\label{e0.5}
\dst\frac12<s<\dst\frac32,\qquad 1<p<\infty.
 \end{eqnarray}
\end{theorem}

Note, that conditions  \eqref{e0.4} and  \eqref{e0.5} are independent of the parameter $p$.

The proof of the foregoing Theorem \ref{t0.3} in \S\, \ref{sect5} is based on the Theorem \ref{t0.3a} and Theorem \ref{t0.4}.
 %
\begin{theorem}\label{t0.3a}
Let $1<p<\infty$, $s>\dst\frac1p$. Let $g_0\in\bW^{s-1/p}_p(\Gamma)$ and $h_0\in\bW^{s-1-1/p}_p(\Gamma)$ be some fixed extensions of the boundary data $g\in\bW^{s-1/p}_p(\Gamma_D)$ and $h\in\bW^{s-1-1/p}_p(\Gamma_N)$ (non-classical formulation), initially defined on the parts of the boundary $\Gamma=\Gamma_D\cup\Gamma_N$.

A solution to the BVP \eqref{e0.1} is represented by the formula
\begin{eqnarray}\label{e0.5a}
u(\cx)=\N_\cC f(\cx)+\W_\Gamma(g_0+\varphi_0)(\cx)-\V_\Gamma(h_0
     +\psi_0)(\cx), \qquad \cx\in\cC.
 \end{eqnarray}
Here $\N_\cC$, $\W_\Gamma$ and $\V_\Gamma$ are the Newton's, double and single layer potentials, defined below $($see \eqref{e1.5}$)$ and $\varphi_0$, $ \psi_0$ in \eqref{e0.5a} are solutions to the following system of pseudodifferential equations
\begin{eqnarray}\label{e0.6}
\begin{array}{l}
\left\{\begin{array}{ll}\dst\frac12\varphi_0-r_N\W_{\Gamma,0}\varphi_0
     +r_N\V_{\Gamma,-1}\psi_0=G_0&\text{on}\quad\Gamma_N,\\[3mm]
\dst\frac12\psi_0+r_D\W^*_{\Gamma,0} \psi_0-r_D\V_{\Gamma,+1}\varphi_0=H_0\qquad
     &\text{on}\quad\Gamma_D, \end{array}\right.
\end{array}\\[2mm]
\label{e0.7}
\begin{array}{c}
\varphi_0\in\wt{\bW}^{s-1/p}_p(\Gamma_N),\quad \psi_0\in\wt{\bW}^{s-1-1/p}_p(\Gamma_D),\\[3mm]
     G_0\in\bW^{s-1/p}_p(\Gamma_N),\qquad H_0\in\bW^{s-1-1/p}_p(\Gamma_D),
\end{array}
\end{eqnarray}
where $G_0$ and $H_0$ are given functions and the participating pseudodifferential operators are defined \eqref{e1.13} in \S\, 1 below.

Vice versa: if $u$  is a solution to the BVP  \eqref{e0.1},  $g:=r_Du^+$, $h:=r_N(\pa_\nub u)^+$ and
$g_0\in\bW^{s-1/p}_p(\Gamma_N)$, $h_0\in\bW^{s-1-1/p}_p(\Gamma_\D)$ are some fixed extensions of $g$ anf $h$ to $\Gamma$, then $\vf_0:=\ r_{\Gamma_D}(u^+ - g_0)$, $\psi_0:=r_{\Gamma_N}((\pa_\nub u)^+ - h_0)$
are solutions to the system \eqref{e0.6}.

The system of boundary pseudodifferential equations \eqref{e0.6} has a unique pair of solutions $\varphi_0\in\bW^{1/2}(\Gamma_N)$ and $\psi_0\in\bW^{-1/2}(\Gamma_D)$ in the classical setting $p=2$, $s=1$.
\end{theorem}

The proof of Theorem \ref{t0.3a} is exposed in \S\, \ref{sect1}.

For the system \eqref{e0.6} we can remove the constraint $s>\dst\frac1p$ and prove the following result for arbitrary $r\in\bR$.
 %
\begin{theorem}\label{t0.4}
Let $1<p<\infty$, $r>-1$.

The system of boundary pseudodifferential equations \eqref{e0.6} is Fredholm in the Sobolev-Slobode\v{c}kii space setting
\begin{subequations}
\begin{eqnarray}\label{e0.6a}
\begin{array}{c}
 \varphi_0\in\wt{\bW}^r_p(\Gamma_N),\quad  \psi_0\in\wt{\bW}^{r-1}_p(\Gamma_D),\\[3mm]
     G_0\in\bW^r_p(\Gamma_N),\qquad H_0\in\bW^{r-1}_p(\Gamma_D)
\end{array}
\end{eqnarray}
and also in the Bessel potential space setting
\begin{eqnarray}\label{e0.6b}
\begin{array}{c}
h, h_0,\varphi_0\in\wt{\bH}^r_p(\Gamma_N),\quad g, g_0,  \psi_0\in\wt{\bH}^{r-1}_p(\Gamma_D),\\[3mm]
     G_0\in\bH^r_p(\Gamma_N),\qquad H_0\in\bH^{r-1}_p(\Gamma_D)
\end{array}
\end{eqnarray}
\end{subequations}
if the following condition holds:
 \begin{eqnarray}\label{e0.8}
p\not=2\quad or \quad p=2\quad and \quad r\not=0,1,2,\ldots.
 \end{eqnarray}
In particular, the system \eqref{e0.6} has a unique solution in both settings \eqref{e0.6a} and \eqref{e0.6b} if:
 \begin{eqnarray}\label{e0.9}
 1<p<\infty, \qquad -1<r<0.
 \end{eqnarray}
\end{theorem}

The proof of the foregoing Theorem \ref{t0.4} in \S\, \ref{sect5} is based on the auxiliary Theorem \ref{t0.5}. To formulate the theorem consider the following model system of boundary integral equations (BIEs)
\begin{eqnarray}\label{e0.10}
&&\hskip-8mm\left\{\begin{array}{ll}\vf(t) + \K^1_{-1}\psi(t)=G(t),\\[3mm]
\psi(t) + \K^1_{-1}\vf(t)=H(t),\qquad  &t\in\bR^+,\end{array}\right.\\[2mm]
&& \varphi,\psi\in\wt{\bW}^{s-1-1/p}_p(\bR^+),\qquad G,\; H\in\bW^{s-1-1/p}_p(\bR^+), \nonumber
\end{eqnarray}
where
\begin{eqnarray}\label{e0.11}
\K^1_{-c}v(t):=\dst\frac1\pi\int_0^\infty\frac{v(\tau)d\tau}{t+c\tau},
     \qquad-\pi<\arg\,c<\pi, \quad v\in\bL_p(\bR^+)
\end{eqnarray}
is a Mellin convolution operator with the kernel homogeneous of order $-1$ (see \cite{Du79,Du84b,Du86,Du82}).
 %
\begin{theorem}\label{t0.5}
Let $1<p<\infty$, $r>-1$. \\
\indent
The system of boundary pseudodifferential equations \eqref{e0.6} is Fredholm in the  Sobolev-Slobode\v{c}kii  \eqref{e0.6a} and Bessel potential \eqref{e0.6b} space settings if the system of boundary integral equations \eqref{e0.10} is locally invertible at $0$ in the Sobolev-Slobode\v{c}kii
\begin{eqnarray}\label{e0.11a}
 \varphi, \psi\in\wt{\bW}^{r-1}_p(\bR^+),\qquad G, H\in\bW^{r-1}_p(\bR^+)
\end{eqnarray}
and the Bessel potential space
\begin{eqnarray}\label{e0.11b}
 \varphi, \psi\in\wt{\bH}^{r-1}_p(\bR^+),\qquad G, H\in\bH^{r-1}_p(\bR^+)
\end{eqnarray}
settings, respectively.
\end{theorem}
 %
\begin{remark}\label{r0.7}
Theorem \ref{t0.5} is proved at the end of \S\, 1. For the proof we apply a quasi-localization of the BVP \eqref{e0.1} with some model BVPs on the half space (see Lemma \ref{l1.4} and Lemma \ref{l1.5}). The constraint $r>-1$ is due to this approach, since the boundary value problems are involved.

In a forthcoming paper will be proved directly the local quasi-equivalence of the equation \eqref{e0.6} and the system \eqref{e0.10} at the points where the Dirichlet and Neumann boundary conditions collide and some simpler equations, which are uniquely solvable, at all other points. Then the constraint $r>-1$ can be dropped and replaced by $r\in\bR$.

Correspondingly, Theorem \ref{t0.4} is also valid for all $r\in\bR$ and the condition \eqref{e0.8} acquires the form
 \[
p\not=2\quad or \quad p=2\quad and \quad r\not=0,\pm1,\pm2,\ldots.
 \]
\end{remark}

A quasi-localization means "freezing coefficients" and "rectifying" underling contours and surfaces. For details of a quasi-localization we refer the reader to the papers \cite{Si65}  and  \cite{CDS03}, where the quasi-localization is well described for singular integral operators and for BVPs, respectively. We also refer to \cite[\S\, 3]{Du15}, where is exposed a short introduction to quasi-localization.

In the present case under consideration we get 3 different model problems by localizing the mixed BVP \eqref{e0.1} to:
\begin{itemize}
  \item[1]
Inner points of $\cC$.
  \item[2]
 Inner points on the boundary $\Gamma_D$ and $\Gamma_N$.
  \item[3]
Points of the boundary $\Gamma$ where different boundary conditions collide (endpoints of $\Gamma_N$ and $\Gamma_D$).
\end{itemize}

The model BVPs obtained by a quasi-localization, are well investigated in the first two cases and such model problems have unique solutions without additional constraints. In the third case we get a mixed BVP on the half plane for the Laplace equation (cf. \eqref{e0.6x} below). The system \eqref{e0.10} is related to this model mixed problem \eqref{e0.6x} just as BVP \eqref{e0.1} is related to the system \eqref{e0.6} (cf.  Lemma \ref{l1.5} below).

The investigation of the boundary integral equation system \eqref{e0.10} is based on recent results on Mellin convolution equations with meromorphic kernels in Bessel potential spaces (see R. Duduchava \cite{Du15}, R. Duduchava and V. Didenko \cite{DD14}).

The symbol $\cB_0^s(\omega)$ of the system \eqref{e0.10} is a continuous function on some infinite rectangle $\mathfrak{R}$ and is responsible for the Fredholm property and the index of the system. This provides necessary and sufficient conditions for the Fredholm property of \eqref{e0.10} which is then used to prove the solvability of the original BVP in the non-classical setting.

 A rigorous analysis of solvability of the above and similar problems with Dirichlet, Neumann, mixed and impedance boundary condition for the Helmholtz and other other elliptic equations are very helpful for a general understanding of elliptic boundary value problems in conical domains (see \cite{KS03,KMR01,No58}).

In \cite{ENS11,ENS14} the authors suggest  another approach to the investigation of the model mixed problem for the Helmholtz equation: they write explicit formulae for a solution with two different methods. But the setting is classical only (the case  $p=2$) and the approach can not be applied to the non-classical setting. Other known results are either very limited to special situations such as the rectangular case \cite{CST04,CST06,MPST93} or apply rather sophisticated analytical methods \cite{KMM05,ZM00}, or are missing a precise setting of appropriate functional spaces (see, e.g., \cite{Ma58,Uf03}).  For the historical survey and for further references we recommend \cite{CK13,ZM00,Va00}.

There is another approach, which can also be applied is the limiting absorption principle, which is based on variational formulation and Lax-Milgram Lemma and its generalizations. Such approach is presented, e.g., in \cite{BT01,BCC12a,BCC12b}. But again, these results are for the classical setting.

In 1960’s there was suggested to solve canonical diffraction problems in Sobolev spaces, based on the recent development in pseudodifferential equations in domains with corners and, more generally,with a Lipschitz boundary. It was popularized by E. Meister \cite{Me85,Me87}, E. Meister and F.-O. Speck \cite{MS79}, W.L. Wendland \cite{WSH79}, A. Ferreira dos Santos \cite{ST89} and their collaborators in the 1980’s. Also see the book of Vasil’ev \cite{Va00} with a considerable list of references. The results are also restricted to the classical setting.

\section{Potential operators and boundary integral equations}
\label{sect1}
\setcounter{equation}{0}

Let $\cS$ be a closed, sufficiently smooth orientable surface in $\mathbb{R}^n$. We use the notation $\mathbb{X}_p^s(\mathcal{S})$ for either the Bessel potential $\mathbb{H}^s_p(\mathcal{S})$ or the Sobolev-Slobode\v{c}kii $\mathbb{W}{}^s_p(\mathcal{S})$ spaces for $\mathcal{S}$ closed or open and a similar notation $\widetilde{\mathbb{X}}_p^s(\mathcal{S})$ for $\mathcal{S}$ open.

Consider the space
 \begin{eqnarray}\label{e1.1}
\mathbb{X}^s_{p,\#}(\mathcal{S}):=\left\{{\varphi}\in\mathbb{X}^s_p
     (\mathcal{S})\;:\;\mbox{\bf{(}}{\varphi},1\mbox{\bf{)}}=0\right\},
 \end{eqnarray}
where  $\mbox{\bf{(}}\cdot,\cdot\mbox{\bf{)}}$ denotes the duality pairing between the adjoint spaces. It is obvious, that $\mathbb{X}^s_{p,\#}(\mathcal{S})$ does not contain constants: if $c_0={\rm const}\in\mathbb{X}^s_{p,\#}(\mathcal{S})$ than
 \[
0=\mbox{\bf{(}} c_0,1\mbox{\bf{)}}=c_0\mbox{\bf{(}}1,1\mbox{\bf{)}}=c_0{\rm mes}\,\mathcal{S}
 \]
and $c_0=0$. Moreover, $\mathbb{X}^s_p(\mathcal{S})$ decomposes into the direct sum
 \begin{eqnarray}\label{e1.2}
\mathbb{X}^s_p(\mathcal{S})=\mathbb{X}^s_{p,\#}(\mathcal{S})+\{{\rm const}\}
 \end{eqnarray}
and the dual (adjoint) space is
 \begin{eqnarray}\label{e1.3}
(\mathbb{X}^s_{p,\#}(\mathcal{S}))^*=\mathbb{X}^{-s}_{p',\#}(\mathcal{S}), \qquad p':=\frac p{p-1}.
 \end{eqnarray}

The following is a part of Theorem 10 proved in \cite{DTT14}.
 %
 \begin{theorem}\label{t1.1a}
Let $\mathcal{S}$ be $\ell$-smooth $\ell=1,2,\ldots$, $1<p<\infty$ and $|s|\leqslant\ell$. Let $\mathbb{X}^s_{p,\#}(\mathcal{S})$  be the same as in \eqref{e1.1}-\eqref{e1.3}.

The Laplace-Beltrami operator $\Delta_\cS:={\bf\rm div}_\mathcal{S}\nabla_\mathcal{S}$ is invertible between the spaces with detached constants
 \begin{eqnarray}\label{e1.4}
\Delta_\cS\;:\;\mathbb{X}^{s+1}_{p,\#}(\mathcal{S})\to\mathbb{X}^{s-1}_{p,\#}(\mathcal{S}),
 \end{eqnarray}
i.e.,  has the fundamental solution $\cK_\cS$ in the setting \eqref{e1.4}.
 \end{theorem}

Let $\cC\subset\cS$ be a subsurface with a smooth boundary $\Gamma:=\partial\cC$.  With the fundamental solution $\cK_\cS$ of the Laplace-Beltrami operator at hand we can consider the standard Newton, single and double layer potentials on the surface $\cC$:
\begin{eqnarray}\label{e1.5}
\begin{array}{l}
\N_\cC v(x):=\dst\int_{\cC}\cK_\cS (x,y)v(y)\,d\sigma\,\\[3mm]
\V_\Gamma v(x):=\dst\int_{\Gamma}\cK_\cS (x,\tau)v(\tau)d\tau,\\[3mm]
\W_\Gamma v(x):=\dst\int_{\Gamma}\pa_{\nub_\Gamma(\tau)}\cK_\cS (x,\tau)v(\tau)d\tau,
     \qquad x\in\cC.
\end{array}
\end{eqnarray}

The potential operators, defined above, have standard boundedness properties
 \begin{eqnarray*}
\N_\cC &:&\bH_{p,\#}^s(\cC)\longrightarrow\bH^{s+2}_{p,\#}(\cC)\, ,\\
\V_\Gamma &:&\bH_{p,\#}^s(\Gamma)\longrightarrow\bH^{s+1+\frac1p}_{p,\#}(\cC)\, ,\\
\W_\Gamma &:&\bH_{p,\#}^s(\Gamma)\longrightarrow\bH^{s+\frac1p}_{p,\#}(\cC)
 \end{eqnarray*}
and any solution to the mixed BVP \eqref{e0.1} in the space $\bH^1_\#(\cC)$ is represented as follows:
\begin{eqnarray}\label{e1.6}
u(x)=\N_\cC f(x)+\W_\Gamma u^+(x)-\V_\Gamma[\pa_{\nub_\Gamma} u]^+(x) \qquad u\in\bH^1_\#(\cC),
     \quad x\in\cC
 \end{eqnarray}
(see \cite{DNS95,Du01}). Densities in \eqref{e1.6} represent the Dirichlet $u^+$ and the Neumann $[\pa_{\nub_\Gamma} u]^+$ traces of the solution $u$ on the boundary.

Since $\bX_p^s=\bX_{p,\#}^s+\{{\rm const}\}$, we can extend layer potentials to the entire space as follows:
 \begin{eqnarray}\label{e1.6b}
 \begin{array}{c}
  \text{for}\quad \varphi=\varphi_0+c, \qquad \varphi_0\in\bX_{p,\#}^s,
  \quad c={\rm const},\\[3mm]
  \text{we set}\quad \V_\Gamma\vf=\V_\Gamma\vf_0+c, \quad \W_\Gamma\vf=\W_\Gamma\vf_0+c, \quad \N_\cC\vf=\N_\cC\vf_0+c,
 \end{array}
   \end{eqnarray}
i.e., by setting $\V_\Gamma c=\W_\Gamma c=c\N_\cC c=c$.
 %
 \begin{lemma}\label{l1.1}
The representation formula \eqref{e1.6} remains valid for a solution in the space $\bH^1(\cC)$, provided the potentials are extended as in \eqref{e1.6b}.
 \end{lemma}
{\bf Proof:} Indeed, since $u=u_0+c$, $u_0\in\bH_{p,\#}^s(\cC)$, $u\in\bH_p^s(\cC)$, we apply the representation formula \eqref{e1.6} for a solution in the space $\bH^1_\#(\cC)$, formula \eqref{e1.6b}, and get the representation formula \eqref{e1.6} for a solution in the space $\bH^1(\cC)$:
\begin{eqnarray}\label{e1.6c}
  \begin{array}{rcl}
  u(x)&\hskip-3mm=&\hskip-3mm u_0(x)+c=\N_\cC f(x)+\W_\Gamma u_0^+(x)-\V_\Gamma[\pa_{\nub_\Gamma} u_0]^+(x)+c\\[3mm]
  &\hskip-3mm=&\hskip-3mm \N_\cC f(x)+\W_\Gamma (u-c)^+(x)-\V_\Gamma[\pa_{\nub_\Gamma} (u-c)]^+(x)+c\\[3mm]
  &\hskip-3mm=&\hskip-3mm \N_\cC f(x)+\W_\Gamma u^+(x)-\V_\Gamma[\pa_{\nub_\Gamma} u]^+(x), \qquad u\in\bH^1(\cC),
  \quad x\in\cC.
  \end{array}
\end{eqnarray}
\vskip-9mm \QED\vskip7mm

\noindent
{\bf Proof of Theorem \ref{t0.3a}:} Let us recall the Plemelji formulae
 \begin{eqnarray}\label{e1.7}
\begin{array}{l}
(\W_\Gamma v)^\pm(t)=\pm\dst\frac12v(t)+\W_{\Gamma,0}v(t),\quad
    (\pa_{\nub_\Gamma}\W_\Gamma\psi)^\pm(t)=\V_{\Gamma,+1}v(t), \\[3mm]
(\pa_{\nub_\Gamma}\V_\Gamma v)^\pm(t)=\mp\dst\frac12v(t)+\W^*_{\Gamma,0}v(t),\quad
(\V_\Gamma v)^\pm(t)=\V_{\Gamma,-1}v(t) ,
 \end{array}
 \end{eqnarray}
where $t\in\pa\Omega_\alpha$ and
\begin{eqnarray}\label{e1.7a}
\begin{array}{rcl}
\V_{\Gamma,-1} v(t) &\hskip-3mm:=&\hskip-3mm\dst\int_\Gamma\cK_\cS (t,\tau) v(\tau)d\tau,\\[3mm]
\W_{\Gamma,0} v(t)&\hskip-3mm:=&\hskip-3mm\dst\int_\Gamma(\pa_{\nub_\Gamma(\tau)}\cK_\cS )(t,\tau)
      v(\tau)d\tau,\\
\W^*_{\Gamma,0}w(t)&\hskip-3mm:=&\hskip-3mm\dst\int_\Gamma(\pa_{\nub_\Gamma(t)}\cK_\cS )(t,\tau)
     w(\tau)d\tau,\\
\V_{\Gamma,+1}w(t)&\hskip-3mm:=&\hskip-3mm\dst\int_\Gamma(\pa_{\nub_\Gamma(t)}\pa_{\nub_\Gamma(\tau)}
     \cK_\cS )(t,\tau)w(\tau)d\tau,\qquad\qquad t\in\Gamma,
 \end{array}
 \end{eqnarray}
are pseudodifferential operators on $\Gamma$, have orders $-1$, $0$, $0$ and $+1$, respectively, and represent the direct values of the corresponding potentials $\V_\Gamma$, $\W_\Gamma$,  $\pa_{\nub_\Gamma}\V_\Gamma$ and $\pa_{\nub_\Gamma}\W_\Gamma$.

Let $g_0\in\bW^{s-1/p}_p(\Gamma)$ and $h_0\in\bW^{s-1-1/p}_p(\Gamma)$ be some fixed extensions of the boundary conditions $g\in\bW^{s-1/p}_p(\Gamma_D)$ and $h\in\bW^{s-1-1/p}_p(\Gamma_N)$ (non-classical formulation), initially defined on the parts of the boundary $\Gamma=\Gamma_D\cup\Gamma_N$. Since the difference between such two extensions belong to the spaces $\wt{\bW}^{s-1/p}_p(\Gamma_N)$ and $\wt{\bW}^{s-1-1/p}_p(\Gamma_D)$ respectively, let us look for two unknown functions $\varphi_0\in\wt{\bW}^{s-1/p}_p(\Gamma_N)$ and $\psi_0\in\wt{\bW}^{s-1-1/p}_p(\Gamma_D)$, such that for $g_0+\vf_0$  and $h_0+\psi_0$ the boundary conditions in \eqref{e0.1} hold on the entire boundary
\begin{eqnarray}\label{e1.8}
\begin{array}{c}
u^+(t)=g_0(t)+\varphi_0(t)=\left\{\begin{array}{ll} g(t)\quad &{\rm if}\quad
     t\in\Gamma_D,\\[3mm]
g_0(t)+\varphi_0(t)\quad &{\rm if}\quad
     t\in\Gamma_N,
\end{array}\right.\\ \\
(\pa_{\nub_\Gamma} u)^+(t)=h_0(t)+\psi_0(t)=\left\{\begin{array}{ll} h_0(t)+\psi_0(t)
     \quad &{\rm if}\quad t\in\Gamma_D,\\[3mm]
h(t)\quad &{\rm if}\quad t\in\Gamma_N,
\end{array}\right.
\end{array}
\end{eqnarray}
provided f $u(x)$ is a solution to the BVP \eqref{e0.1}.

By introducing the boundary values of a solution \eqref{e1.8} to the BVP \eqref{e0.1} into the representation formula \eqref{e1.6c} (see Lemma \ref{l1.1}) we get the following representation of a solution:
\begin{eqnarray}\label{e1.9}
u(x)=\N_\cC f(x)+\W_\Gamma[g_0+\varphi_0](x)-\V_\Gamma[h_0+\psi_0](x),\qquad x\in\cC,
\end{eqnarray}
where
 \[
\hskip-5mm g_0\in\bW^{s-1/p}_p(\Gamma), \; h_0\in\bW^{s-1-1/p}_p(\Gamma),\;
\varphi_0\in\wt{\bW}^{s-1/p}_p(\Gamma_N), \; \psi_0\in\wt{\bW}^{s-1-1/p}_p(\Gamma_D).
 \]

By applying the Plemelji formulae \eqref{e1.7} to \eqref{e1.9} and taking into account \eqref{e1.8} we get the following:
 \[
\begin{array}{r}
\left\{\begin{array}{l}
g_0(t)+\varphi_0(t)=u^+(t)={(\N_\cC f)^+}+\dst\frac12(g_0(t)+\varphi_0(t))\\
\hskip20mm+\W_{\Gamma,0}[g_0+\varphi_0](t)-\V_{\Gamma,-1}[h_0+\psi_0](t),\\[2mm]
h_0(t)+\psi_0(t)=(\pa_{\nub_\Gamma} u)^+(t)={(\pa_{\nub_\Gamma}\N_\cC f)^+}+\V_{\Gamma,+1}[g_0+\varphi_0](t)\\
\hskip20mm+\dst\frac12(h_0(t) +\psi_0(t))-\W^*_{\Gamma,0}[h_0+\psi_0](t),\qquad t\in\Gamma.
\end{array}\right.
\end{array}
 \]
If we apply the restriction operator $r_D$ to $\Gamma_D$ to the first equation in the obtained system and the restriction operator $r_N$ to $\Gamma_N$ to the second one, we obtain the system \eqref{e0.6}, where
\begin{eqnarray}\label{e1.13}
\begin{array}{c}
G_0:=r_N\left[(\N_\cC f)^+-\dst\frac12g_0+\W_{\Gamma,0}g_0-\V_{\Gamma,-1}
     h_0\right]\in\bW^{s-1/p}_p(\Gamma_N),\\[2mm]
H_0:=r_D\left[(\pa_{\nub_\Gamma}\N_\cC f)^+-\dst\frac12h_0+\V_{\Gamma,+1}
     g_0-\W^*_{\Gamma,0}h_0\right]\in\bW^{s-1-1/p}_p(\Gamma_D).
\end{array}
\end{eqnarray}

Thus, we have proved the inverse assertion of Theorem \ref{t0.3a}: if $u$ is a solution to the BVP \eqref{e0.1}, the functions $\vf_0$ and $\psi_0$ are solutions to the system \eqref{e0.6}.

The direct assertion is easy to prove:
 \begin{itemize}
 \item
The function in \eqref{e1.6c} represented by the potentials, satisfies the equation \eqref{e0.1}.
\item
If $\vf_0$ and $\psi_0$ are solutions to the system \eqref{e0.6}, using Plemelji formulae \eqref{e1.7}  it can easily be verified that $u$ in \eqref{e1.6c} satisfies the boundary conditions in \eqref{e0.1}.
\end{itemize}

The existence and the uniqueness of a solution to the BVP \eqref{e0.1} in the classical setting \eqref{e0.2} is stated in Theorem \ref{t0.3}, while for the system \eqref{e0.6} it follows from the equivalence with the BVP \eqref{e0.1}.     \QED

The remainder of the paper is devoted to the proof of solvability properties of the system \eqref{e0.6} in the non-classical setting \eqref{e0.3}.

Consider the following equation on the 2-dimensional Euclidean space
\begin{equation}\label{e0.6m}
\Delta u=f^0\qquad\text{on}\quad \bR^2, \qquad u\in\mathbb{H}^{s}_p(\bR^2),\quad f^0\in\mathbb{H}^{s-2}_p(\bR^2),
\end{equation}
also the model Dirichlet
\begin{eqnarray}\label{e0.4m}
\left\{\begin{array}{ll}
\Delta u(x)=f_0(x),\qquad & x\in\bR^2_+, \\[0.2cm]
u^+(t)=g_0(t),     \qquad & t\in\partial\bR^2_+=\bR, \\[0.2cm]
\end{array}\right.
\end{eqnarray}
the model Neumann
\begin{eqnarray}\label{e0.5m}
\left\{\begin{array}{ll}
\Delta u(x)=f_0(x),\qquad & x\in\bR^2_+, \\[0.2cm]
-(\partial_2u)^+(t)=h_0(t), \qquad &t\in\partial\bR^2_+=\bR,
\end{array}\right.
\end{eqnarray}
and the model mixed
\begin{eqnarray}\label{e0.6x}
\left\{\begin{array}{ll}
\Delta u(x)=f_1(x),\qquad & x\in \bR^2_+, \\[0.2cm]
     u^+(t)=g_1(t),     \qquad & t\in\bR^-:=(-\infty,0), \\[0.2cm]
-(\partial_{2}u)^+(t)=h_1(t),\qquad
     &t\in \bR^+:=(0,\infty),
\end{array}\right.
\end{eqnarray}
boundary value problems for the Laplace equation on the upper half plane $\bR^2_+:=\bR\times\bR^+$, where $\pa_{\nub_\Gamma}=-\pa_2$ is the normal derivative on the boundary of $\bR^2_+$.

The BVPs \eqref{e0.4m} and \eqref{e0.5m} will be treated in the non-classical setting:
\begin{equation}\label{e0.7m}
\begin{array}{r}
f_0\in\widetilde{\mathbb{H}}^{s-2}_p(\bR^2_+)\cap\widetilde{
    \mathbb{H}}^{-1}_0(\bR^2_+),\quad g_0\in\mathbb{W}^{s-1/p}_p(\bR), \quad h_0\in\mathbb{W}^{s-1-1/p}_p(\bR),\\[3mm]
1<p<\infty, \qquad s>\dst\frac1p
\end{array}
\end{equation}
and the BVP \eqref{e0.6x} will be treated in the non-classical setting:
\begin{equation}\label{e0.7x}
\begin{array}{r}
f_1\in\widetilde{\mathbb{H}}^{s-2}_p(\bR^2_+)\cap\widetilde{\mathbb{H}}^{-1}_0(\bR^2_+),\quad
     g_1\in\mathbb{W}^{s-1/p}_p(\bR^-),\quad h_1\in\mathbb{W}^{s-1-1/p}_p(\bR^+),\\[3mm]
     1<p<\infty, \qquad s>\dst\frac1p.
\end{array}
\end{equation}

 %
 \begin{proposition}\label{p1.7}
The BVPs \eqref{e0.4m}, \eqref{e0.5m} have unique solutions in the setting \eqref{e0.7m} and the Laplace equation in the setting \eqref{e0.6m} has a unique solution as well.
\end{proposition}
{\bf Proof:} The assertion is a well-known classical result, available in many textbooks on partial differential equations (see e.g. \cite{HW08}).  \QED

As a paticular case of Theorem \ref{t0.1} (can easily be proved with the Lax-Milgram Lemma) we have the following.
 %
\begin{proposition}\label{p1.4}
The mixed BVP \eqref{e0.6x} has a unique solution $u$ in the classical weak setting
\begin{equation*}
u\in\mathbb{H}^1(\bR^2_+),\quad f_1\in\widetilde{\mathbb{H}}^{-1}_0(\bR^2_+),
     \quad g_1\in\mathbb{H}^{1/2}(\bR^+), \quad h_1\in\mathbb{H}^{-1/2}(\bR^-),
\end{equation*}
\end{proposition}
 %
 \begin{lemma}\label{l1.4}
The BVP \eqref{e0.1} is Fredholm  in the non-classical setting \eqref{e0.3} if the model mixed BVP \eqref{e0.6x} is locally Fredholm (ie., is locally invertible) at $0$ in the non-classical setting \eqref{e0.7x}.
\end{lemma}
{\bf Proof:}  We apply quasi-localization of the boundary value problem \eqref{e0.1} in the more general non-classical setting \eqref{e0.3}, which includes the classical setting  \eqref{e0.2} as a particular case (see \cite{CDS03,Du84a}) for details of quasi-localization of boundary value problems and also \cite{DiS08,GK79,Si65} for general results on localization and quasi-localization.

By quasi-localization at the point $\omega\in\ov\cC$ we first localize to the tangential plane $\bR^2(\omega)$ (tangential half plane $\bR^2_+(\omega)$)  to $\cC$ at $\omega\in\cC$ (at $\omega\in\Gamma=\partial\cC$, respectively). The differential operators remain the same
 \begin{equation}\label{e0.x7}
\begin{array}{c}
\Delta_{\bR^2}:=\dst\sum_{j=1}^3\cD_j^2, \quad \cD_j=\pa_j-\nu_j\pa_{\nub},\\ \pa_{\nub}=\dst\sum_{j=1}^3\nu_j\cD_j, \quad \pa_{\nub_\Gamma}=\dst\sum_{j=1}^3\nu_{\Gamma,j}\cD_j,
\end{array}
\end{equation}
but the normal vector $\nub(\omega)$  to the tangent plane $\bR^2$ and the normal vector $\nub_\Gamma(\omega)$  to the boundary of the tangent plane $\bR(\omega)=\partial\bR^2_+(\omega)$ are now constant. Next we rotate the tangent planes $\bR^2(\omega)$ and $\bR^2_+(\omega)$ to match them with the plane $\bR^2$ and $\bR^2_+$. The normal vector fields $\nub(\omega)$ will transform into $\nub=(0,0,1)$ and $\nub_\Gamma(\omega)=(0,-1,0)$. The rotation is an isomorphism of the spaces $\bW^r_p(\bR^2(\omega))\to\bW^r_p(\bR^2)$, $\bW^r_p(\bR^2_+(\omega))\to\bW^r_p(\bR^2_+)$, $\wt\bW^r_p(\bR^2_+(\omega))\to\wt\bW^r_p(\bR^2_+)$ etc. and transforms the operators in \eqref{e0.x7} into the operators
 \begin{equation*}
\begin{array}{c}
\Delta_{\bR^2(\omega)}\;\to\;\Delta:=\dst\sum_{j=1}^2\pa_j^2, \quad \cD_j\;\to\;\pa_j,
     \quad j=1,2,\quad, \cD_3\;\to\;0,\\
\pa_{\nub(\omega)}\;\to\; \pa_3, \quad \pa_{\nub_\Gamma(\omega)}\;\to\;-\pa_2
\end{array}
\end{equation*}
and we get \eqref{e0.6m}, \eqref{e0.4m}, \eqref{e0.5m}, \eqref{e0.6x} as a local representatives of BVP \eqref{e0.1}.

For the  BVP \eqref{e0.1} in the non-classical setting \eqref{e0.3} we get the following local quasi-equivalent equations and BVPs at different points of the surface $\omega\in\ov\cC$:
 \begin{itemize}
 \item[i.]
The equation \eqref{e0.6m} at $0$ if $\omega\in\cC$ is an inner points of the surface;
 \item[ii.]
The Dirichlet BVP \eqref{e0.4m} in the non-classical setting \eqref{e0.7m} at $0$  if $\omega\in\Gamma_D$;
 \item[iii.]
The Neumann BVP \eqref{e0.5m} in the non-classical setting \eqref{e0.7m} at $0$ if
$\omega\in\Gamma_N$;
 \item[iv.]
The mixed BVP \eqref{e0.6x} in the non-classical setting \eqref{e0.7x} at $0$ if $\omega\in\ov{\Gamma_D}\cap\ov{\Gamma_N}$ is one of two points of collision of different boundary conditions.
 \end{itemize}

The main conclusion of the present theorem on Fredholm properties of BVPs \eqref{e0.1} and \eqref{e0.6x} follows from Proposition \ref{p1.7} and the general theorem on quasi-localizaion  (see \cite{CDS03,Du84a,DiS08,GK79,Si65}): {\em The BVP \eqref{e0.1}, \eqref{e0.3} is Fredholm if all local representatives \eqref{e0.6m}, \eqref{e0.4m}, \eqref{e0.5m} and \eqref{e0.6x} in non-classical settings are locally Fredholm (i.e., are locally invertible)}.        \QED

Now we concentrate on the model mixed BVP \eqref{e0.6x}.To this end let us recall that the function
 \begin{eqnarray*}
    \cK_\Delta(x):=\frac1{2\pi}\ln|x|
 \end{eqnarray*}
is the fundamental solution to the Laplace's equation in two variables
 \begin{eqnarray}\label{e1.25a}
 \begin{array}{cc}
\Delta\cK_\Delta(x) = \delta(x), \qquad x\in\bR^2,\\[3mm]
\Delta=\partial^2_1+\partial^2_2=\partial^2_\nub+\partial^2_\ell.
\end{array}
 \end{eqnarray}

From \eqref{e1.25a} follows the equality
 \[
\delta=\Delta\cK_\Delta=\partial^2_\nub\cK_\Delta +\partial^2_\ell\cK_\Delta,
 \]
which we use to prove the following:
\begin{eqnarray}\label{e1.25b}
\partial_{\nub(x)}\partial_{\nub(y)}\cK_\Delta(x-y)=-\partial^2_{\nub(y)}\cK_\Delta(x-y)
     =-\delta(x-y)+\partial^2_{\ell(y)}\cK_\Delta(x-y).
 \end{eqnarray}

Applying the latter equality \eqref{e1.25b}, we represent the hypersingular operator $\V_{\bR,+1}$ as follows
 \begin{eqnarray}\label{e1.25c}
&&\hskip-10mm\V_{\bR,+1}\vf(t):=\dst\int_{\bR}\pa_{\nub(t)}\pa_{\nub(\tau)}\cK_\Delta(t-\tau)
     \vf(\tau)d\tau=-\vf(t)+\dst\int_{\bR}\pa^2_{\ell(\tau)}\cK_\Delta(t-\tau)\vf(\tau)d\tau\nonumber\\
&&\hskip-0mm=-\vf(t)-\dst\int_{\bR}\pa_\tau\cK_\Delta(t-\tau)\pa_\tau\vf(\tau)d\tau, \quad t\in\bR,
\end{eqnarray}
since $\pa_{\ell(\tau)}=\pa_\tau$ on $\bR$ and for the tangential differential operator $\pa_\ell$ on arbitrary smooth contour $\Gamma$ the following "partial integration" formula is valid (see \cite{Du01,DMM06}):
 \[
\dst\int_\Gamma\pa_{\ell(\tau)}\psi(\tau)\vf(\tau)d\sigma=-\dst\int_{\Gamma}
     \psi(\tau)\pa_{\ell(\tau)}\vf(\tau)d\sigma.
 \]

We can define standard layer potential operators, the Newton, the single and the double layer potentials respectively (cf. \eqref{e1.5})
 \begin{eqnarray}\label{e1.20x}
\N_{\bR^2_+}v(x)&\hskip-3mm:=&\hskip-3mm\dst\frac1{2\pi}\dst\int_{\bR^2_+}
     \ln|x-y| v(y)\,dy,\nonumber\\[3mm]
\V_{\bR} v(x)&\hskip-3mm:=&\hskip-3mm\dst\frac1{2\pi}\dst\int_{\bR}\ln|x-\tau|
      v(\tau)d\tau,\nonumber\\
\W_{\bR} v(x)&\hskip-3mm:=&\hskip-3mm-\dst\frac1{2\pi}\dst\int_{\bR}\pa_{2}
     \ln|(x_1,x_2)-(\tau,y_2)|\Big|_{y_2=0} v(\tau)d\tau\\
&\hskip-3mm=&\hskip-3mm-\dst\frac1{2\pi}\dst\int_\bR\pa_{2}
     \ln\sqrt{(x_1-\tau)^2+(x_2-y_2)^2}\Big|_{y_2=0} v(\tau)d\tau\nonumber\\[3mm]
&\hskip-3mm=&\hskip-3mm\dst\frac1{2\pi}\dst\int_\bR\frac{x_2 v(\tau)d\tau}{(x_1-\tau)^2
     +x_2^2}, \qquad x=(x_1,x_2)^\top\in\bR^2_+.\nonumber
\end{eqnarray}

The pseudodifferential operators on $\V_{\bR,-1}$, $\W_{\bR,0}$, $\W^*_{\bR,0}$ and $\V_{\bR,+1}$, associated with the layer potentials (see \eqref{e1.7a}), acquire the form
\begin{eqnarray}\label{e1.27}
\begin{array}{rcl}
\V_{\bR,-1}v(t)&\hskip-3mm:=&\hskip-3mm\dst\frac1{2\pi}\dst\int_\bR\ln|t-\tau|
     v(\tau)d\tau,\\[3mm]
\W_{\bR,0} v(t)&\hskip-3mm:=&\hskip-3mm\lim\limits_{x_2\to0}\dst\frac1{2\pi}\dst\int_\bR
     \frac{x_2 v(\tau)d\tau}{(x_1-\tau)^2+x_2^2}=0,\qquad \W^*_{\bR,0}v(t)=0,\\[3mm]
 \end{array}
 \end{eqnarray}

By using the representation \eqref{e1.25c} we find the following:
\begin{eqnarray*}
\begin{array}{rcl}
\V_{\bR,+1} v(t)
 &\hskip-3mm=&\hskip-3mm-v(t)-\dst\frac1{2\pi}\dst\int_{\bR}\pa_\tau\ln|t-\tau|\pa_\tau v(\tau)d\tau\\[3mm]
 &\hskip-3mm=&\hskip-3mm-v(t)+\dst\frac1{2\pi}\dst\int_{\bR}\dst\frac{t-\tau}{(t-\tau)^2}\pa_\tau v(\tau)d\tau\\[3mm]
  &\hskip-3mm=&\hskip-3mm-v(t)+\dst\frac1{2\pi}\dst\int_{\bR}\dst\frac{\pa_\tau v(\tau)d\tau}{t-\tau}, \quad t\in\bR
 \end{array}
 \end{eqnarray*}
and the Plemelji formulae \eqref{e1.7} acquire the form
 \[
\begin{array}{c}
(\W_{\bR} v)^\pm(t)=\pm\dst\frac12 v(t),\qquad
     -(\pa_{y_2}\V_\bR v)^\pm)(t)=\mp\dst\frac12v(t),\\[3mm]
-(\pa_{y_2}\W_{\bR}v)^\pm(t)=\V_{\bR,+1}v(t)\qquad
(\V_{\bR} v)^\pm(t)=\V_{\bR,-1} v(t) \qquad t\in\bR.
 \end{array}
 \]

Now we prove the following.
 %
\begin{lemma}\label{l1.5}
Let $1<p<\infty$, $s>\dst\frac1p$. Let $g^0_1\in\bW^{s-1/p}_p(\bR)$ and $h^0_1\in\bW^{s-1-1/p}_p(\bR)$  be some fixed extensions of the boundary conditions $g_1\in\bW^{s-1/p}_p(\bR^-)$ and $h_1\in\bW^{s-1-1/p}_p(\bR^+)$ $($non-classical formulation \eqref{e0.7x}$)$, initially defined on the parts of the boundary $\bR=\bR^-\cup\bR^+$.

A solution to the BVP \eqref{e0.6x} is represented by the formula
\begin{eqnarray}\label{e1.28}
u(x)=\N_{\bR^2_+}f(x)+\W_\bR(g^0_1+\varphi^0)(x)-\V_\bR(h^0_1+\psi^0)(x),\qquad x\in\bR^2
 \end{eqnarray}
(cf. \eqref{e1.20x} for the potential operators) and $\varphi^0$ and $\psi^0$ are solutions to the system of pseudodifferential equations
\begin{eqnarray}\label{e1.22}
&&\begin{array}{l}
\left\{\begin{array}{ll}\dst\frac12\varphi^0 -  r_+\W_{\bR,0}\varphi^0 + r_+\V_{\bR,-1}\psi^0=G_1
     \qquad \text{on}\quad \bR^+\\[3mm]
\dst\frac12\psi^0 + r_-\W^*_{\bR,0}\psi^0 - r_-\V_{\bR,+1}\varphi^0=H_1 \qquad \text{on}\quad \bR^-, \end{array}\right.
\end{array}  \\[2mm]
\label{e1.23}
&&\varphi^0, \R_*\psi^0\in\wt\bW^{s-1-1/p}_p(\bR^+),\quad G_1, \R_* H_1\in\bW^{s-1-1/p}_p(\bR^+),
\end{eqnarray}
where $r_+$ and $r_-$ are the restriction operators from the axes $\bR$ to the semi- axis $\bR^+$ and $\bR^-$.

The system of boundary pseudodifferential equations \eqref{e1.22} has a unique pair of solutions $\varphi^0$ and $\psi^0$ in the classical setting $p=2$, $s=1$.
\end{lemma}
 {\bf Proof:} By repeating word by word the proof of Theorem \ref{t0.3a}, we prove the equivalence via the representation formulae \eqref{e1.28} of the BVP \eqref{e0.6x} in the non-classical setting \eqref{e0.7x} and of the system \eqref{e1.22}.

 The existence and uniqueness of a solution to the BVP \eqref{e0.6x} in the classical setting \eqref{e0.7x} is stated in Proposition \ref{p1.4}, while for the system \eqref{e1.22} it follows from the proved equivalence with the BVP \eqref{e0.6x}. \QED
 %
\begin{lemma}\label{l1.5a}
Let $1<p<\infty$, $s>\dst\frac1p$.

The system of boundary pseudodifferential equations \eqref{e1.22}
is locally invertible at $0$ if and only if the system \eqref{e0.10} is locally invertible at $0$ in the non-classical setting  \eqref{e0.11a} and the space parameters are related as follows: $r=s-\frac1p>0$.
\end{lemma}
 {\bf Proof:} Due to the equalities \eqref{e1.27} $r_+\W_{\bR,0}\varphi^0=0$,  $r_-\W^*_{\bR,0}\psi^0=0$ and the equation in \eqref{e1.22} acquires  the form
\begin{eqnarray*}
\left\{\begin{array}{ll}\dst\frac12\varphi^0(t) + \dst\frac1{2\pi}\dst\int_{\bR^-}\ln|t-\tau|\psi^0(\tau)d\tau=G_1(t),
     \qquad t\in\bR^+,\\[3mm]
\dst\frac12\psi^0(t)  -\dst\frac1{2\pi}\dst\int_{\bR^+}\frac{(\pa_\tau\varphi^0)(\tau)
     d\tau}{t-\tau}=H_1(t), \qquad  t\in\bR^-. \end{array}\right.
\end{eqnarray*}
Multiply both equations by 2, apply to the first equation the differentiation $\pa_t$, replace $\vf:=\pa_t\vf^0$,  apply to the second equation the reflection $\R_* v(t)=v(-t)$ and replace $\psi=\R_*\psi^0$, also under the integral. We get the following
\begin{eqnarray*}
\left\{\begin{array}{ll}\varphi(t) + \dst\frac1\pi\dst\int_{\bR^+}\pa_t\ln(t+\tau)\psi(\tau)d\tau
     =\varphi(t) + \dst\frac1\pi\dst\int_{\bR^+}\dst\frac{\psi(\tau)d\tau}{t+\tau} =2\pa_tG_1(t)=:G(t), \\[5mm]
\psi(t)  +\dst\frac1\pi\dst\int_{\bR^+}\frac{\varphi(\tau)d\tau}{t+\tau}=2H_1(-t)=:H(t), \qquad  t\in\bR^+ \end{array}\right.
\end{eqnarray*}
and the obtained equation coincides with the system \eqref{e0.10}.

To prove the local equivalence at $0$ of the systems  \eqref{e1.22} and \eqref{e0.10} note, that the multiplication by $2$ and the reflection
 \[
\R_*\;\;:\bW^r_p(\bR^+)\to\bW^r_p(\bR^-),\qquad
\R_*\;:\;\wt\bW^r_p(\bR^+)\to\wt\bW^r_p(\bR^-)
 \]
are invertible operators since $\R_*^2=I$ and $\R_*^{-1}=\R_*$ and, therefore, are locally invertible at $0$.

The differentiation
 \[
\pa_t:=\dst\frac d{dt}\;:\;\bW^r_p(\bR^+)\to\bW^{r-1}_p(\bR^+),\qquad
\pa_t\;:\;\wt\bW^r_p(\bR^+)\to \wt\bW^{r-1}_p(\bR^+)
 \]
is locally invertible at any finite point $x\in\bR$ because the operators
 \[
\pa_t-iI\;:\;\bW^r_p(\bR^+)\to\bW^{r-1}_p(\bR^+),\qquad
\pa_t+iI\;:\;\wt\bW^r_p(\bR^+)\to \wt\bW^{r-1}_p(\bR^+)
 \]
are isomorphisms (represent Bessel potentials, see Theorem \ref{t4.8} below, \cite[Lemma 5.1]{Du79} and \cite{Es81}). On the other hand, the embeddings
 \[
iI\;:\;\bW^r_p(\bR^+)\to\bW^{r-1}_p(\bR^+),\qquad
iI\;:\;\wt\bW^r_p(\bR^+)\to \wt\bW^{r-1}_p(\bR^+)
 \]
are locally compact due to the Sobolev's embedding theorem and the compact perturbation does not influences the local invertibility.  \QED

\noindent
{\bf Proof of Theorem \ref{t0.5}:}  By Theorem \ref{t0.3a} the system \eqref{e0.6} is Fredholm  in the  Sobolev-Slobode\v{c}kii  space setting \eqref{e0.6a} if the BVP \eqref{e0.1} is Fredholm in the non-classical setting \eqref{e0.3}. On the other hand, by Lemma  \ref{l1.4} the BVP \eqref{e0.1} is Fredholm in the non-classical setting \eqref{e0.3} if the BVP \eqref{e0.6x} is locally invertible at $0$ in the non-classical setting \eqref{e0.7x}. And, finally, by Lemma  \ref{l1.5} and Lemma  \ref{l1.5a} the BVP \eqref{e0.6x} is locally invertible in the non-classical setting \eqref{e0.7x} if the system of boundary integral equations \eqref{e0.10} is locally invertible at $0$ in the Sobolev-Slobode\v{c}kii space setting \eqref{e0.11a}. This accomplishes the proof of the first part of the assertion, concerning the solvability in the Sobolev-Slobode\v{c}kii  space settings \eqref{e0.6a} and \eqref{e0.11a}.

The second part of the assertion, concerning the solvability in the Bessel potential space settings  \eqref{e0.6b} and \eqref{e0.11b}, follows from the first part and Proposition \ref{p4.3}, exposed below and proved in \cite{Du15,DD14}, which states that these solvability properties are equivalent.     \QED

\section{Fourier convolution operators in the Bessel potential spaces $\bH^s_p(\bR^+)$}
\label{sect2}
\setcounter{equation}{0}

To formulate the next theorem we need to introduce Fourier convolution and Bessel potential operators.

For the spaces of scalar, vector and matrix functions we will use the same notation if this will not lead to a confusion. For example, $\mathbb{L}_{\infty,loc}(\mathbb{R})$ might be the space of locally bounded functions either scalar, but also vector or matrix valued functions; this will be clear from the context.

Let $a\in\mathbb{L}_{\infty,loc}(\mathbb{R})$ be a locally bounded $m\times m$ matrix function. The Fourier convolution operator (FCO) with the symbol $a$ is defined by
\[
    W^0_a:=\mathcal{F}^{-1}a\mathcal{F}.
 \]
Here
 \begin{equation*}
\cF u(\xi):=\dst\int_{\bR^n}e^{i\xi x}u(x)dx,\quad \xi\in\bR^n,
 \end{equation*}
is the Fourier transform and
 \begin{equation*}
\cF^{-1}v(\xi):=\dst\frac{1}{(2\pi)^n}\int_{\bR^n}e^{-i\xi x}v(\xi)d\xi, \quad x\in\bR^n,
 \end{equation*}
is its inverse transform. If the operator
\begin{equation*}
    W^0_a:\mathbb{H}^s_p(\mathbb{R})\longrightarrow \mathbb{H}^{s-r}_p(\mathbb{R})
\end{equation*}
is bounded, we say that $a$ is an $\mathbb{L}_p$-multiplier of order $r$ and use "$\mathbb{L}_p$-multiplier" if the order is $0$. The set of all $\mathbb{L}_p$-multipliers of order $r$ (of order $0$) is denoted by $\mathfrak{M}^r_p(\mathbb{R})$ (by $\mathfrak{M}_p(\mathbb{R})$, respectively).
Let
 \[
\wt{\mathfrak{M}}^r_p(\mathbb{R}):=\bigcap_{p-\ve<q<p+\ve}\mathfrak{M}^r_q(\mathbb{R}),
\qquad\wt{\mathfrak{M}}_p(\mathbb{R}):=\bigcap_{p-\ve<q<p+\ve}\mathfrak{M}_q(\mathbb{R}).
 \]
Note, that $\wt{\mathfrak{M}}^r_p(\mathbb{R})$ and $\wt{\mathfrak{M}}_p(\mathbb{R})$ are independent of $\ve$ because, due to interpolation theorem $\mathfrak{M}^r_{p_0}(\mathbb{R})\subset\mathfrak{M}^r_{p_-}(\mathbb{R})\bigcap
\mathfrak{M}^r_{p_+}(\mathbb{R})$ for all $1<p_-<p_0<p_+<\infty$.

For an $\mathbb{L}_p$-multiplier of order $r$, $a\in\mathfrak{M}^r_p(\mathbb{R})$, the Fourier convolution operator (FCO) on the semi-axis $\mathbb{R}^+$ is defined by the equality
\begin{equation}\label{e4.7}
W_a=r_+W^0_a\; :\; \widetilde{\mathbb{H}}^s_p(\mathbb{R}^+)\longrightarrow
     \mathbb{H}^{s-r}_p(\mathbb{R}^+)
\end{equation}
where $r_+:=r_{\mathbb{R}^+}:\mathbb{H}^s_p(\mathbb{R})\longrightarrow \mathbb{H}^s_p(\mathbb{R}^+)$ is the restriction operator to the semi-axes
$\mathbb{R}^+$.

We did not use the parameter $s\in\mathbb{R}$ in the definition of the class of multipliers $\mathfrak{M}^r_p(\mathbb{R})$ . This is due to the fact that $\mathfrak{M}^r_p(\mathbb{R})$ is independent of $s$: if the operator $W_a$ in \eqref{e4.7} is bounded for some $s\in\mathbb{R}$, it is bounded for all other values of $s$. Another definition of the multiplier class $\mathfrak{M}^r_p(\mathbb{R})$ is written as follows: $a\in\mathfrak{M}^r_p(\mathbb{R})$ if and only if $\lambda^{-r}a\in\mathfrak{M}_p(\mathbb{R})=\mathfrak{M}^0_p(\mathbb{R})$, where $\lambda^r(\xi):=(1+|\xi|^2)^{r/2}$. This assertion is one of the consequences of Theorem \ref{t4.8} below.

Consider  the Bessel potential operators defined as follows
 \begin{equation}\label{e4.1}
  \begin{array}{l}
\mathbf{\Lambda}_\gamma^r=W^0_{\lambda^r_\gamma}\;:\;\widetilde{\mathbb{H}
     }^s_p(\mathbb{R}^+)\rightarrow\widetilde{\mathbb{H}}^{s-r}_p(\mathbb{R}^+),
     \,\\[3mm]
\mathbf{\Lambda}_{-\gamma}^r=r_+W^0_{\lambda^r_{-\gamma}}\ell\;:\;\mathbb{H}^s_p(
     \mathbb{R}^+)\rightarrow\mathbb{H}^{s-r}_p(\mathbb{R}^+)\, ,\\[3mm]
\lambda^r_{\pm\gamma}(\xi):=(\xi\pm\gamma)^r,\qquad\xi\in\mathbb{R}, \qquad {\rm Im}\,\gamma>0
 \end{array}
 \end{equation}
for a non-negative $s\geqslant0$. Here $\ell\;:\;\mathbb{H}^s_p(\mathbb{R}^+)\to\mathbb{H}^s_p(\mathbb{R})$ is some extension operator. In \eqref{e4.7} there is no need of any extension operator since the space $\widetilde{\mathbb{H}}^s_p(\mathbb{R}^+)$ is automatically embedded in $\bH^s_p(\bR)$ provided functions are extended by $0$.

For a negative $s<0$ the Bessel potential operators $\mathbf{\Lambda}_{\pm\gamma}^r$ are defined by the duality between the spaces.
 %
\begin{theorem}\label{t4.8}
Let  $1<p<\infty$. Then:
\begin{enumerate}
\vskip+0.15cm
\item
For any $r,s\in\mathbb{R}$, $\gamma\in\mathbb{C}$, ${\rm Im}\,\gamma>0$ the Bessel potential
operators \eqref{e4.1} arrange isomorphisms of the corresponding spaces $($see {\rm \cite{Du79, Es81}}$)$ and are independent of the choice of an extension operator $\ell:\mathbb{H}^s_p(\mathbb{R}^+)\longrightarrow\mathbb{H}^s_p (\mathbb{R})$.

\vskip+0.15cm
\item
For any  operator  $\mathbf{A}:\widetilde{\mathbb{H}}^s_p(\mathbb{R}^+) \longrightarrow \mathbb{H}^{s-r}_p(\mathbb{R}^+)$ of order $r$, the following diagram is commutative
\begin{equation}\label{e4.10}
\begin{array}{ccc}\widetilde{\mathbb{H}}^s_p(\mathbb{R}^+) & \stackrel{\mathbf{A}}{\longrightarrow}
    &\mathbb{H}^{s-r}_p(\mathbb{R}^+)\\
\uparrow\mathbf{\Lambda}^{-s}_\gamma & &\downarrow \mathbf{\Lambda}_{-\gamma}^{s-r}\\
\mathbb{L}_p(\mathbb{R}^+) & \stackrel{\mathbf{\Lambda}_{-\gamma}^{s-r}\mathbf{A}
    \mathbf{\Lambda}^{-s}_\gamma}{\longrightarrow}& \mathbb{L}_p(\mathbb{R}^+). \end{array}
 \end{equation}

The diagram \eqref{e4.10} provides an equivalent lifting of the
operator $\mathbf{A}$ of order $r$ to the operator
$\mathbf{\Lambda}_{-\gamma}^{s-r}\mathbf{A}\mathbf{\Lambda}^{-s}_\gamma:\mathbb{L}_p(
\mathbb{R}^+)\longrightarrow\mathbb{L}_p(\mathbb{R}^+)$ of order~$0$.
\vskip+0.15cm
\item
For any bounded convolution operator $W_a:\mathbb{H}^s_p(\mathbb{R}^+) \longrightarrow \mathbb{H}^{s-r}_p(\mathbb{R}^+)$ of order $r$ and for any pair of complex numbers $\gamma_1, \gamma_2$ such that ${\rm Im}\,\gamma_j>0$, $j=1,2$, the lifted operator
\begin{equation}\label{e4.11}
 \begin{array}{c}
\mathbf{\Lambda}_{-\gamma_1}^\mu W_a\mathbf{\Lambda}_{\gamma_2}^\nu
     =W_{a_{\mu,\nu}}\;:\;\mathbb{H}^{s+\nu}_p(\mathbb{R}^+) \longrightarrow\mathbb{H}^{s-r-\mu}_p(\mathbb{R}^+),\\[2mm]
a_{\mu,\nu}(\xi):=(\xi-\gamma_1)^\mu a(\xi)(\xi+\gamma_2)^\nu
 \end{array}
\end{equation}
is again a Fourier convolution.

In particular,  the lifted operator  $W_{a_0}$ in $\mathbb{L}_p$-spaces, $\mathbf{ \Lambda}_{-\gamma}^{s-r}W_a\mathbf{\Lambda}^{-s}_\gamma:\mathbb{L}_p(\mathbb{R}^+)
\longrightarrow\mathbb{L}_p(\mathbb{R}^+)$ has the symbol
\begin{equation*}
a_{s-r,-s}(\xi)=\lambda^{s-r}_{-\gamma}(\xi)a(\xi)\lambda^{-s}_\gamma(\xi)
    =\Big(\frac{\xi-\gamma}{\xi+\gamma}\Big)^{s-r}\,\frac{a(\xi)}{(\xi+i)^r}\,.
\end{equation*}
\end{enumerate}
\end{theorem}
 %
\begin{remark}\label{r4.12}
For any pair of multipliers $a\in\mathfrak{M}^r_p(\mathbb{R})$, $b\in\mathfrak{M}^s_p(\mathbb{R})$
the corresponding convolution operators on the full axes $W^0_a$ and $W_b^0$ have the property $W^0_aW^0_b=W^0_bW^0_a=W^0_{ab}$.

For the corresponding Wiener-Hopf operators on the half axis a similar equality
\begin{equation}\label{e4.13}
W_aW_b=W_{ab}
\end{equation}
is valid if at least one of the following conditions hold: the function $a(\xi)$ has an analytic extension in the lower half plane or the function $b(\xi)$ has an analytic extension in the upper half plane (see \cite{Du79}).

Note, that actually \eqref{e4.11} is a consequence of \eqref{e4.13}.
\end{remark}

Let $\dR:=\bR\cup\{\infty\}$ denote the one point compactification of the real axis $\bR$ and $\overline{\bR}:=\bR\cup\{\pm\infty\}$-the two point compactification of $\bR$. By $C(\dR)$ (by $C(\overline{\bR})$, respectively) we denote the space of continuous functions $g(x)$ on $\bR$ which have the same limits at the infinity $g(-\infty)=g(+\infty)$ (limits at the infinity might differ $g(-\infty)\not=g(+\infty)$). By $PC(\dR)$ is denoted the space of piecewise-continuous functions on $\dR$, having limits $a(t\pm0)$ at all points $t\in\dR$, including infinity.
 %
\begin{proposition}[Lemma 7.1, \cite{Du79} and Proposition 1.2, \cite{Du87}]\label{p2.5} Let $1<p<\infty$, $a\in C(\dR{}^+)$, $b\in C(\dR)\cap\wt{\mathfrak{M}}_p(\dR)$ and $a(\infty)= b(\infty)=0$. Then the operators $aW_b,W_b\,aI:\mathbb{L}_p(\mathbb{R}^+)\longrightarrow \mathbb{L}_p(\mathbb{R}^+)$ are compact.

Moreover, these operators are compact in all Bessel potential and Besov spaces, where they are bounded, due to the Krasnoselskij interpolation theorem for compact operators.
\end{proposition}
 %
\begin{proposition}[Lemma 7.4, \cite{Du79} and Lemma 1.2, \cite{Du87}]\label{p2.7}
Let $1<p<\infty$ and let $a$ and $b$ satisfy at least one of the following conditions:
\begin{itemize}
\vskip+0.15cm
\item[(i)] $a\in C(\overline{\mathbb{R}}{}^+)$, $b\in\wt{\mathfrak{M}}_p(\mathbb{R})\cap PC(\overline{\mathbb{R}})$,
\vskip+0.15cm
\item[(ii)] $a\in PC(\overline{\mathbb{R}}{}^+)$, $b\in C\wt{\mathfrak{M}}_p(\overline{\mathbb{R}})$.
\end{itemize}
Then the commutants $[aI,W_b]$ are compact operators in the space $\mathbb{L}_p(\mathbb{R}^+)$ and also, due to Krasnoselskij interpolation theorem for compact operators, in all Bessel potential and Besov spaces, where they are bounded.
\end{proposition}

\section{Mellin convolution operators in the space $\bH^s_p(\bR^+)$}
\label{sect3}
\setcounter{equation}{0}

In this section we expose auxiliary results from \cite{Du15} (also see \cite{Du79,Du87,DD14}), which are essential for the investigation of boundary integral equations from the foregoing section.

Let $a(\xi)$  be a $N\times N$ matrix function $a\in C\mathfrak{M}^0_p(\bR)$, continuous on the real axis $\mathbb{R}$ with the only possible jump at infinity.  Consider a Mellin convolution operator $\mathfrak{M}^0_a$ with the symbol $a$ in the Bessel potential spaces
 \begin{eqnarray*}
\mathfrak{M}^0_a:=\cM^{-1}_\beta a\cM_\beta\;:\;\widetilde{\mathbb{H}}^s_p(
     \mathbb{R}^+)\longrightarrow\mathbb{H}^s_p(\mathbb{R}^+),\quad s\in\mathbb{R},
 \end{eqnarray*}
where
 \begin{eqnarray*}
 \begin{array}{c}
\mathcal{M}_\beta v(\xi):=\displaystyle\int_0^\infty \tau^{\beta-i\xi}v(\tau)\frac{d\tau}\tau,\quad
     \xi\in\mathbb{R},\\[1ex]
\mathcal{M}^{-1}_\beta u(t):=\displaystyle\frac1{2\pi}
     \displaystyle\int_{-\infty}^{\infty}t^{i\xi-\beta}
     u(\xi)d\xi,\qquad t\in\mathbb{R}^+,
 \end{array}
 \end{eqnarray*}
are the Mellin transformation and the inverse to it.

The most important example of a Mellin convolution operator is an integral operator of the form
 \begin{eqnarray}\label{e2.2}
\mathfrak{M}^0_a\mathbf{u}(t):=c_0\mathbf{u}(t)+\frac{c_1}{\pi i}
    \int_0^\infty\frac{\mathbf{u}(\tau)\,d\tau}{\tau-t}+\int_{0}^\infty
    \mathcal{K}\left(\frac t\tau\right)\mathbf{u}(\tau)\frac{d\tau}{\tau}
 \end{eqnarray}
with $n\times n$ matrix coefficients and $n\times n$ matrix kernel
 \begin{eqnarray}\label{e2.4}
\int_0^\infty t^{\bt-1}\cK(t)dt<\infty, \quad 0<\bt<1.
  \end{eqnarray}
Then $\mathfrak{M}^0_a$ is a bounded operator in the weighted Lebesgue space of vector functions
 \begin{eqnarray}\label{e2.3}
\fM^0_a\;:\;\bL_p(t^\gamma,\mathbb{R}^+)\longrightarrow\bL_p(t^\gamma,\mathbb{R}^+),\\[3mm] \bt:=\frac{1+\gamma}p, \quad 1<p<\infty,\quad -1<\gamma<p-1,\nonumber
 \end{eqnarray}
endowed with the norm
 \[
\|u|\bL_p(t^\gamma,\mathbb{R}^+)\|:=\left[\int_0^\infty t^\gamma|u(t)|^pdt
     \right]^{1/p}
 \]
(cf. \cite{Du79}). The symbol of the operator \eqref{e2.2} is the Mellin transform of the kernel
 \begin{eqnarray*}
a_\bt(\xi)&\hskip-3mm:=&\hskip-3mm c_0+c_1\coth\,\pi\left(i\beta
     +\xi\right)+\cM_\bt\cK(\xi)\nonumber\\
&\hskip-3mm:=&\hskip-3mm c_0+c_1\coth\,\pi\left(i\beta+\xi\right)
     +\int_{0}^{\infty}t^{\bt-i\xi}\cK(t)\frac{dt}t,\quad \xi\in\bR.
 \end{eqnarray*}

Obviously,$ \mathfrak{M}^0_a\mathfrak{M}^0_b\vf =\mathfrak{M}^0_{ab}\vf$ for $\vf\in C^\infty_0(\mathbb{R}^+)$.
 %
 \begin{theorem}\label{t2.1}
Let $1<p<\infty$ and $-1<\gamma<p-1$ (or $0\leqslant p\leqslant\infty$ provided $c_1=0$ in \eqref{e2.2}). The following three properties are equivalent:
\begin{itemize}
  \item [i.]
  Operator $\fM_a^0$ in \eqref{e2.2}--\eqref{e2.3}  is Fredholm;
  \item [ii.]
   The symbol of the operator is invertible (is elliptic)
 \begin{eqnarray*}
\inf_{\xi\in\bR}\left|\det\,a_\beta(\xi)\right|>0;
 \end{eqnarray*}
  \item [iii.]
The operator is invertible and the inverse operator is $\fM_{a^{-1}}^0$.
\end{itemize}
\end{theorem}
 %
\begin{proposition}[Lemma 7.4, \cite{Du79} and Lemma 1.2, \cite{Du87}]\label{p2.7a}
Let $1<p<\infty$ and let $a$ and $b$ satisfy at least one of the following conditions:
\begin{itemize}
\vskip+0.15cm
\item[(i)] $a\in C(\overline{\mathbb{R}}{}^+)$, $b\in\wt{\mathfrak{M}}_p(\mathbb{R})\cap PC(\overline{\mathbb{R}})$,
\vskip+0.15cm
\item[(ii)] $a\in PC(\overline{\mathbb{R}}{}^+)$, $b\in C\wt{\mathfrak{M}}_p(\overline{\mathbb{R}})$.
\end{itemize}
Then the commutants $[aI,\mathfrak{M}^0_b]$ are compact operators in the space $\mathbb{L}_p(\mathbb{R}^+)$ and also, due to Krasnoselskij interpolation theorem for compact operators, in all Bessel potential and Besov spaces, where they are bounded.
\end{proposition}

Things are different in the Bessel potential spaces \ if\  compared\  with the\  Lebesgue\  spaces. Let us recall some results from \cite[\S\, 2]{Du15}. Consider meromorphic functions in the complex plane $\mathbb{C}$, vanishing at infinity
 \begin{eqnarray}\label{e2.6a}
 \begin{array}{r}
\mathcal{K}(t):=\dst\sum_{j=0}^N\dst\frac{d_j}{(t-c_j)^{m_j}}
 \end{array}
 \end{eqnarray}
with poles at $c_0,c_1,\ldots\in\mathbb{C}\setminus\{0\}$, complex coefficients $d_j\in\mathbb{C}$ and $m_j\in\mathbb{N}$.
 \begin{definition}[see \cite{Du15}]\label{d2.2}
We call a kernel $\mathcal{K}(t)$ in \eqref{e2.6a} admissible if for those poles $c_0,\ldots,c_\ell$ which belong to the positive semi-axes $\arg\,c_0=\cdots =\arg\,c_\ell=0$, the corresponding multiplicities are one, i.e., $m_0=\cdots=m_\ell=1$.
 \end{definition}

For example: The Mellin convolution operator
\begin{eqnarray*}
\K^m_cv(t):=\dst\frac1\pi\int_0^\infty\frac{\tau^{m-1}v(\tau)d\tau}{(t-c\tau)^m},
     \qquad  0<\arg\,c<2\pi,\quad t\in\bR^+, \quad v\in\bL_p(\bR^+)
 \end{eqnarray*}
has an admissible kernel for arbitrary $m=1,2,\ldots$ if $m=1$ as soon as $c$ is real $\arg\,c=0$.
 %
\begin{proposition}[see \cite{Du15}, Corollary 2.3, Theorem 2.4]\label{p2.4}
Let $1<p<\infty$ and $-1<\gamma<p-1$ (or $1\leqslant p\leqslant\infty$ provided $c_1=0$ in \eqref{e2.2}) and $\mathcal{K}(t)$ in \eqref{e2.6a} be an admissible kernel. Then the Mellin convolution
\[
\mathfrak{M}^0_{a_\beta}\mathbf{u}(t):=c_0\mathbf{u}(t)+\int_{0}^\infty
    \mathcal{K}\left(\frac t\tau\right)\mathbf{u}(\tau)\frac{d\tau}{\tau}
\]
is a bounded operator in the Lebesgue space $\mathbb{L}_p(\mathbb{R}^+,t^\gamma)\to \mathbb{L}_p(\mathbb{R}^+,t^\gamma)$ and, also, in the Bessel potential spaces $\mathfrak{M}^0_{a_\beta}\;:\;\widetilde{\mathbb{H}}^s_p(\mathbb{R}^+)\to
\mathbb{H}^s_p(\mathbb{R}^+)$ for all $s\in\bR$.
 \end{proposition}

The next theorem provides the lifting of the Mellin convolution operator from a pair of Bessel potential spaces to the Lebesgue spaces.
 %
\begin{theorem}[\cite{Du15}, Theorem 4.1]\label{t4.7}
Let  $0<\arg\,c<2\pi$, $0<\arg \gamma<\pi$ and  $r,s\in\mathbb{R}$, $1<p<\infty$. Then the operator $K^1_c\;:\;\widetilde{\mathbb{H}}^s_p (\mathbb{R}^+)\to\mathbb{H}^s_p(\mathbb{R}^+)$ is lifted equivalently to the operator
 \begin{eqnarray*}
\A^{1,s}_c:=\mathbf{\Lambda}^s_{-\gamma}\K^1_c\mathbf{\Lambda}^{-s}_{\gamma}
     \;:\;\mathbb{L}_p(\mathbb{R}^+)\to\mathbb{L}_p(\mathbb{R}^+),
 \end{eqnarray*}
where
 \[
\A^{1,s}_c= c^{-s}\K^1_c W_{g^s_{-c\gamma,\gamma}}, \qquad c^{-s}:=|c|^{-s}e^{-\arg\,c\,ri}
 \]
if only $0<\arg(-c\gamma)<\pi$.

If $0<\arg\,(c\,\gamma)<\pi$, choose any $\gamma_0\in\bC$ such that $0<\arg \gamma_0<\pi$ and $0<\arg(-c\,\gamma_0)<\pi$ (such a choice of $\gamma_0$ is possible since $c$ is not a real constant $\arg\,c\not=0$). Then
 \begin{eqnarray*}
\A^{1,s}_c=c^{-s}W_{g^s_{-\gamma,-\gamma_0}\cdot}\K^1_c
     W_{g^s_{-c\gamma_0,\gamma}}
     =\mathbf{K}^1_cW_{g^s_{-\gamma,-\gamma_0}g^s_{-c\gamma_0,\gamma}}
     \!+\!\mathbf{T},\\
g^s_{-c\gamma_0,\gamma}(\xi):=\left(\displaystyle\frac{\xi-c\gamma_0}{\xi
     +\gamma}\right)^s,\quad g^s_{-\gamma,-\gamma_0}(\xi):=\left(\displaystyle\frac{\xi-\gamma}{\xi
     -\gamma_0}\right)^s,
 \end{eqnarray*}
where $\mathbf{T}\;:\;\mathbb{L}_p(\mathbb{R}^+)\to\mathbb{L}_p(\mathbb{R}^+)$ is a compact operator.
\end{theorem}

\section{Investigation of a lifted Mellin convolution operator}
\label{sect4}
\setcounter{equation}{0}

The results of the foregoing two sections together with results on a Banach algebra generated by Mellin and Fourier convolution operators (see \cite{Du87}) allow the investigation of lifted Mellin convolution operators. For this we need to write the symbol of a model operator
\begin{equation}\label{e4.1a}
\mathbf{A}:=d_0I+\sum_{j=1}^nd_j\mathbf{K}^1_{c_j}\;:\;\wt{\mathbb{H}}^s_p(\mathbb{R}^+)\to
     \mathbb{H}^s_p( \mathbb{R}^+),
\end{equation}
where $\mathbf{K}^1_{c_1},\ldots,\mathbf{K}^1_{c_n}$ are admissible Mellin convolution operators.

To expose the symbol of the operator \eqref{e4.1a}, consider the infinite clockwise oriented ``rectangle'' $\mathfrak{R}:=\Gamma_1\cup\Gamma_2^-\cup\Gamma_2^+ \cup\Gamma_3$, where (cf. Figure~1)
$$  \Gamma_1:=\overline{\mathbb{R}}\times\{+\infty\},\;\;\Gamma^\pm_2:=\{\pm\infty\}\times\overline{\mathbb{R}}^+,\;\;
            \Gamma_3:=\overline{\mathbb{R}}\times\{0\}.         $$

\setlength{\unitlength}{0.4mm}
\vskip7mm
\hskip15mm
\begin{picture}(300,140)
\put(-00,40){\epsfig{file=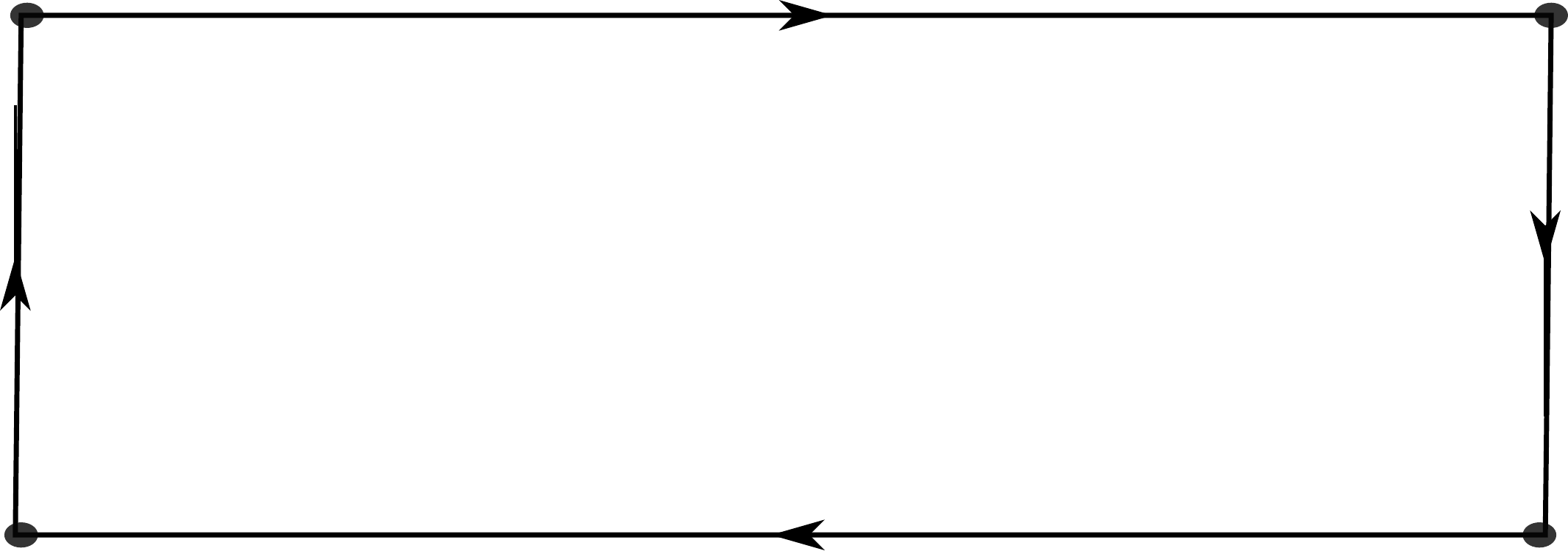,height=40mm, width=80mm}}
\put(80,50){\makebox(0,0)[lc]{$(0,\xi)$}}
\put(80,130){\makebox(0,0)[lc]{$(\infty,\xi)$}}
\put(80,33){\makebox(0,0)[lc]{$\Gamma_3$}}
\put(80,145){\makebox(0,0)[lc]{$\Gamma_1$}}
\put(-13,90){\makebox(0,0)[lc]{$\Gamma^-_2$}}
\put(8,90){\makebox(0,0)[lc]{$(\eta,-\infty)$}}
\put(203,90){\makebox(0,0)[lc]{$\Gamma^+_2$}}
\put(160,90){\makebox(0,0)[lc]{$(\eta,+\infty)$}}
\put(-10,145){\makebox(0,0)[lc]{$(\infty,-\infty)$}}
\put(170,35){\makebox(0,0)[lc]{$(0,+\infty)$}}
\put(-10,35){\makebox(0,0)[lc]{$(0,-\infty)$}}
\put(170,145){\makebox(0,0)[lc]{$(\infty,+\infty)$}}
\put(0,15){\makebox(0,0)[lc]{The domain $\mathfrak{R}$ of definition of the symbol $\mathcal{A}^s_p(\omega)$.}}
\end{picture}

\noindent
According to \cite{DD14} the symbol $\mathcal{A}^s_p(\omega)$ of the operator $\mathbf{A}$ is
\begin{equation}\label{e4.2}
\mathcal{A}^s_p(\omega):=d_0\mathcal{I}^s_p(\omega)+\sum_{j=1}^nd_j\mathcal{K}^{1,s}_{c_j,p}(\omega),
\end{equation}
where
\begin{subequations}
 \begin{eqnarray}\label{e4.3a}
\mathcal{I}^s_p(\omega)&\hskip-3mm:=&\hskip-3mm\begin{cases}
    g^s_{-\gamma,\gamma,p}(\infty,\xi), & \omega=(\infty,\xi)\in\overline{\Gamma}_1,
    \\[1ex]
\left(\displaystyle\frac{\eta-\gamma}{\eta+\gamma}\right)^{\mp s}, &
     \omega=(\eta,\pm\infty)\in\Gamma^\pm_2, \\[1ex]
     e^{\pi si}, &\omega=(0,\xi)\in\overline{\Gamma}_3, \qquad \xi,\eta\in\mathbb{R},\end{cases}\\[2ex]
 &&\hskip-20mm g^s_{-\gamma,\gamma,p}(\infty,\xi):=\frac{e^{2\pi si}+1}2
     +\frac{e^{2\pi si}-1}{2i}\cot\pi\Big(\frac1p-i\xi\Big)=e^{\pi si}\frac{\sin\pi\Big(\frac1p+s-i\xi\Big)}
     {\sin\pi\Big(\frac1p-i\xi\Big)}, \quad \xi\in\mathbb{R},\nonumber
 \end{eqnarray}
 \begin{eqnarray}\label{e4.3b}
\mathcal{K}^{1,s}_{c,p}(\omega):=\begin{cases}
    \displaystyle\frac{e^{-i\pi(\frac1p-i\xi-1)}c^{\frac1p-i\xi-s-1}}{
    \sin\pi(\frac1p-i\xi)},&\omega=(\infty,\xi)
    \in\overline{\Gamma}_1,\\[1ex]
    0, &\omega=(\eta,\pm\infty)\in\Gamma^\pm_2,\\[1ex]
    \displaystyle\frac{e^{-i\pi(\frac1p-i\xi-1)}c^{\frac1p-i\xi-s-1}}{
    \sin\pi(\frac1p-i\xi)},&\omega=(0,\xi)
    \in\overline{\Gamma}_3,\end{cases}\\[1.5ex]
0<\arg\,c<2\pi,\quad 0<|\arg(c\,\gamma)|<\pi,\quad 0<\arg\gamma<\pi
     \nonumber
 \end{eqnarray}
and $c^\delta=|c|^\delta e^{i\delta\arg\,c}$, $\delta\in\bR$.
\end{subequations}

Note, that the Mellin convolution operator $\mathbf{K}^1_{-1}$,
 \[
\begin{array}{c}
\displaystyle
\mathbf{K}^1_{-1}\varphi(t)=\mathbf{K}^1_{e^{i\pi}}\varphi(t)=\frac1\pi
     \int\limits_0^\infty\displaystyle\frac{\varphi(\tau)\,d\tau}{t+\tau}
     =\mathfrak{M}^0_{k_p}\varphi(t),\qquad k_p(\xi)=\displaystyle\frac{1}{\sin\pi\left(\frac1p-i\xi \right )},
 \end{array}
 \]
which we encounter in applications (see \eqref{e0.10} and  Lemma \ref{l1.5}), has a rather simple symbol in the Bessel potential space $\mathbb{H}^s_p(\mathbb{R}^+)$: From  \eqref{e4.3b} follows that:
\begin{equation}\label{e4.4}
\hskip-5mm\mathcal{K}^{1,s}_{-1,p}(\omega):=\begin{cases}
     \displaystyle\frac{e^{-\pi si}}{\sin\pi(\beta-i\xi)},&\omega =(\infty,\xi)\in\overline{\Gamma}_1,\\
     0, &\omega=(\eta,\pm\infty))\in\Gamma^\pm_2,\\
     \displaystyle\frac{e^{-\pi si}}{\sin\pi(\beta-i\xi)},&\omega =(0,\xi)\in\overline{\Gamma}_3. \end{cases}
 \end{equation}

The image of the function $\det\mathcal{A}^s_p(\omega)$, $\omega\in\mathfrak{R}$ is a closed curve in the complex plane (easy to check analyzing the symbol in \eqref{e4.3a}-\eqref{e4.3b}).  Hence, if the symbol is elliptic, i.e. if
 \[
    \inf_{\omega\in\mathfrak{R}} \big|\det\mathcal{A}^s_p(\omega)\big|>0,
 \]
 the increment of the argument $(1/2\pi)\arg \mathcal{A}^s_p(\omega)$ when $\omega$ ranges through $\mathfrak{R}$ in the direction of orientation, is an integer. It is called the winding number or the index of the curve $\Gamma:=\{z\in \mathbb{C}: z=\det\mathcal{A}_p(\omega),\;\omega \in \mathfrak{R}\}$ and is denoted by ${\rm ind}\,\det \mathcal{A}^s_p$.

Propositions \ref{p4.1}-\ref{p4.3}, exposed below, are well known and will be applied in the next section in the proof of main theorems.
 %
\begin{proposition}[\cite{Du15} and Theorem 5.4, \cite{DD14}]\label{p4.1}
Let $1<p<\infty$, $s\in\mathbb{R}$. The operator
\begin{equation}\label{e4.5}
\mathbf{A}:\widetilde{\mathbb{H}}{}^s_p(\mathbb{R}^+)\longrightarrow
     \mathbb{H}^s_p (\mathbb{R}^+)
\end{equation}
defined in \eqref{e4.1} is Fredholm if and only if its symbol $\mathcal{A}^s_p(\omega)$ defined in  \eqref{e4.2}, \eqref{e4.3a}-- \eqref{e4.3b}, is elliptic. If $\mathbf{A}$ \, is Fredholm, then
 \[
{\rm Ind}\mathbf{A}=-{\rm ind}\det\mathcal{A}^s_p.
 \]

The operator $\mathbf{A}$ in \eqref{e4.5} is locally invertible at $0$
if and only if its symbol $\mathcal{A}^s_p(\omega)$ is elliptic on the set $\Gamma_1$ only: $\inf_{\omega\in\Gamma_1} \big|\det\mathcal{A}^s_p(\omega)\big|>0$.
\end{proposition}
 %
 \begin{proposition}[\cite{Du15,DD14}]\label{p4.2}
Let $1<p<\infty$, $s\in\mathbb{R}$ and let $\mathbf{A}$ be defined by \eqref{e4.1}. If the operator $\mathbf{A}\;:\;\widetilde{\mathbb{H}}{}^s_p (\mathbb{R}^+) \longrightarrow\mathbb{H}^s_p (\mathbb{R}^+)$ is Fredholm (is invertible) for all $a\in(s_0,s_1)$ and $p\in(p_0,p_1)$, where $-\infty<s_0<s_1 <\infty$,  $1<p_o<p_1<\infty$, then $\A$ is Fredholm (is invertible, respectively) in the  Sobolev-Slobode\v{c}kii space setting
 \[
\mathbf{A}\;:\;\widetilde{\mathbb{W}}{}^s_p (\mathbb{R}^+) \longrightarrow
     \mathbb{W}^s_p (\mathbb{R}^+),\qquad \text{for all}\quad s\in(s_0,s_1) \quad\text{and}\quad p\in(p_0,p_1)
 \]
and has the same index
 \[
{\rm Ind}\,\mathbf{A}=-{\rm ind}\,\det\,\mathcal{A}^s_p.
 \]
\end{proposition}
 %
 \begin{proposition}[\cite{Du73,DNS95}]\label{p4.3}
Let two pairs of parameter-dependent Banach spaces $\fB^s_1$ and $\fB^s_2$, $s_1<s<s_2$, have intersections $\fB^{s'}_j\cap\fB^{s''}_j$ dense in $\fB^{s'}_j$ and in $\fB^{s''}_j$ for all $j=1,2$, $s',s''\in(s_1,s_2)$.

If a linear bounded operator $A\;:\;\fB^s_1\to\fB^s_2$ is Fredholm for all $s\in(s_1,s_2)$, it has the same kernel and co-kernel for all values of this parameter $s\in(s_1,s_2)$.

In particular, If $A\;:\;\fB^s_1\to\fB^s_2$ is Fredholm for all $s\in(s_1,s_2)$ and is invertible for only one value $s_0\in(s_1,s_2)$, it is invertible for all values of  this parameter $s\in(s_1,s_2)$.
\end{proposition}

\section{Investigation of the boundary integral equations}
\label{sect5}
\setcounter{equation}{0}

The proof of Theorem \ref{t0.4} (see below) is based, besides Theorem \ref{t0.5}, on the following theorem.
 %
\begin{theorem}\label{t5.1}
Let $1<p<\infty$, $r\in\bR$.

The system of the boundary pseudodifferential equations \eqref{e0.10} is Fredholm in the  Sobolev-Slobode\v{c}kii space setting \eqref{e0.11a} and in the Bessel potential space setting \eqref{e0.11b} if and only if the condition \eqref{e0.8} holds. The system \eqref{e0.10} has a unique solution in both settings \eqref{e0.11a} and \eqref{e0.11b} if the condition \eqref{e0.9} holds.
\end{theorem}
{\bf Proof:} Let us write the equation \eqref{e0.10} in an operator form
 \begin{subequations}
\begin{eqnarray}\label{e5.2a}
&\M\Phi=\F, \quad \M:=\left[\begin{array}{cc} I & \K^1_{-1}\\
     \K^1_{-1}& I \end{array}\right],\\
\label{e5.2b}
&\Phi:=\left(\begin{array}{c}\vf\\ \psi\end{array}\right)\in
     \widetilde{\bW}{}^r_p(\bR^+),\qquad {\bf F}:=\left(\begin{array}{c}G\\
     H\end{array}\right)\in\bW^r_p(\bR^+),\\
\label{e5.2c}
&\Phi:=\left(\begin{array}{c}\vf\\ \psi\end{array}\right)\in
     \widetilde{\bH}{}^r_p(\bR^+),\qquad {\bf F}:=\left(\begin{array}{c}G\\
     H\end{array}\right)\in\bH^r_p(\bR^+)
\end{eqnarray}
\end{subequations}
and apply Proposition \ref{p4.1} to the investigation of equation \eqref{e5.2a} in the setting \eqref{e5.2c}.

Due to formulae \eqref{e4.3a} and \eqref{e4.4} the symbol of $\M$ on $\Gamma_1$ reads
 \begin{eqnarray}\label{e5.3}
\cM^r_p(\omega)=
    \left[\begin{array}{cc}  e^{\pi r i}\dst\frac{\sin\pi(\Xi+r)}{\sin\pi\Xi} &\dst\frac{ e^{-\pi r i}}{\sin\pi\Xi}\\[3mm]
     \dst\frac{ e^{-\pi r i}}{\sin\pi\Xi}&  e^{\pi r i}\dst\frac{\sin\pi(\Xi+r)}{\sin\pi\Xi} \end{array}\right],\quad \omega=(\infty,\xi)\in\overline{\Gamma_1},
 \end{eqnarray}
where  $\Xi:=\dst\frac1p-i\xi$,  $\xi\in\bR$, $\eta\in\bR^+$. We have dropped the information about the symbol $\cM^r_p(\omega)$ on the contours $\Gamma^\pm_2$ and $\Gamma_3$ because, due to Theorem \ref{t0.5}, we are interested only in the local invertibility of the operator $\M$ at $0$. This information, due to the concluding part of the Proposition \ref{p4.1}, is contained in the symbol $\cM^r_p(\omega)$ on the contour $\Gamma_1$ only.

According the formula \eqref{e5.3} the symbol $\cM^r_p(\infty,\xi)$ is elliptic on the contour $\Gamma_1$ if and only if
 \[
\det\cM^r_p(\infty,\xi)=\dst\frac{ e^{2\pi r i}\sin^2\pi\left(\dst\frac1p+r -i\xi\right) -  e^{-2\pi r i}}{\sin^2\pi\left(\dst\frac1p
     -i\xi\right)}\not=0, \qquad \omega\in\Gamma_1
 \]
or, equivalently,
 \[
\sin^2\pi\left(\dst\frac1p+r -i\xi\right)\not= e^{-4\pi r i}=\cos\,4\pi r-i\sin4\pi r\qquad \text{for all}\quad \xi\in\bR.
 \]

The symbol is non-elliptic if
 \[
 \sin4\pi r=0 \quad\text{and}\quad  \sin^2\pi\left(\dst\frac1p+r\right)=\cos\,4\pi r=\pm1.
 \]
The latter equation has the following solutions
 \begin{eqnarray}\label{e5.1}
4\pi r=2\pi k \quad\text{and}\quad \sin^2\pi\left(\dst\frac1p+\frac k2\right)=1, \qquad k=0\pm1,\ldots,
 \end{eqnarray}
because for $4\pi r=2k+1$ the equation  $\sin^2\pi\left(\dst\frac1p+r\right)=-1$ has no solution. Equation \eqref{e5.1} decomposes into the following two equations for even and odd $k$:
 \[
 \begin{array}{c}
r=k, \quad  \sin^2\dst\frac\pi p=1\quad \Rightarrow\quad r=k,\quad p=2, \quad k=0,\pm1,\ldots,\\
r=k+\dst\frac12,\quad \cos^2\dst\frac\pi p=1\quad \Rightarrow \quad r=k+\dst\frac12,\quad p=1, \quad k=0,\pm1,\ldots.
\end{array}
 \]

Due to Proposition \ref{p4.2} the operator $\M$ in \eqref{e5.2a} is Fredholm in the setting \eqref{e5.2b} if and only if the same condition \eqref{e0.8} holds.

From \eqref{e0.8} follows that if conditions \eqref{e0.9} hold, the operator $\M$ is Fredholm in both settings  \eqref{e5.2b} and  \eqref{e5.2c}. On the other hand, for the values $p=2$, $r=-1/2$, which also satisfy the conditions  \eqref{e0.9}, the operator $\M$ is invertible (see the concluding assertion in Lemma \ref{l1.5}). Then, due to Proposition \ref{p4.3}, $\M$ is invertible in both settings  \eqref{e5.2b} and  \eqref{e5.2c} for all those $r$ and $p$ which satisfy \eqref{e0.9}.         \QED

\noindent
{\bf Proof of Theorem \ref{t0.4}:} The Fredholm criterion \eqref{e0.8} for the system of boundary pseudodifferential equations \eqref{e0.6} in the settings \eqref{e0.6a} and \eqref{e0.6b} is a direct consequence of Theorem \ref{t0.5} and Theorem \ref{t5.1}.

From \eqref{e0.8} follows that, if conditions \eqref{e0.9} hold, the operator $M_0$, corresponding to the system \eqref{e0.6}, is Fredholm in both settings  \eqref{e0.6a} and \eqref{e0.6b}. On the other hand, for the values $p=2$, $r=-1/2$, which also satisfy the conditions  \eqref{e0.9}, the operator $\M_0$ is invertible (see the concluding assertion in Theorem \ref{t0.3a}). Then, due to Proposition \ref{p4.3}, $\M_0$ is invertible in both settings  \eqref{e0.6a} and  \eqref{e0.6b} for all those $r$ and $p$ which satisfy \eqref{e0.9}.         \QED

\noindent
{\bf Proof of Theorem \ref{t0.3}:} Due to Theorem \ref{t0.3a} and Theorem \ref{t0.4} the BVP \eqref{e0.1} is Fredholm if the system \eqref{e0.6} in the non-classical setting \eqref{e0.7} is, provided $r=\dst\frac1p-s$, i.e., if the condition \eqref{e0.8} holds with $r=\dst\frac1p-s$ (cf.  the condition \eqref{e0.8}), which is the same condition as \eqref{e0.4}.

From \eqref{e0.4} follows that if conditions \eqref{e0.5} hold, the BVP \eqref{e0.1} is Fredholm in the non-classical setting \eqref{e0.3}. On the other hand, for the values $p=2$, $s=1$, which also satisfy the conditions  \eqref{e0.5}, the BVP \eqref{e0.1} has a unique solution (see Theorem \ref{t0.1}). Then, due to Proposition \ref{p4.3}, the BVP \eqref{e0.1} has a unique solution in the non-classical setting \eqref{e0.3} for all those $s$ and $p$ which satisfy \eqref{e0.5}.         \QED

\noindent
{\bf R. Duduchava}, {\em A.Razmadze Mathematical Institute,  Tbilisi State University, Tamarashvili str. 6, Tbilisi 0177, Georgia.} \ {\sf email: RolDud@gmail.com}\\
\\
{\bf M. Tsaava}, {\em A.Razmadze Mathematical Institute,  Tbilisi State University, Tamarashvili str. 6, Tbilisi 0177, Georgia.}\ \ {\sf email: m.caava@yahoo.com}
\end{document}